\documentclass[11pt]{article}
\usepackage{mathrsfs}
\usepackage{amssymb}

\def\hang{\hangindent\parindent}
\def\textindent#1{\indent\llap{#1\enspace}\ignorespaces}
\def\re{\par\hang\textindent}

\title{Looking for Gr\"obner Basis Theory for (Almost)\\ Skew 2-Nomial Algebras
\thanks{Project supported by
the National Natural Science Foundation of China (10571038,
10971044).\newline e-mail: huishipp@yahoo.com}}

\author{Huishi Li\\
{\small Department of Applied Mathematics}\\
{\small College of Information Science and Technology}\\
{\small Hainan University}\\
{\small  Haikou 570228, China}}

\date{}

\begin{document}
\maketitle
\begin{center}
\begin{minipage}{120mm}
{\small {\bf Abstract.}  In this paper, we introduce (almost) skew
2-nomial algebras and look for a one-sided or two-sided  Gr\"obner
basis theory for such algebras at a modest level. That is, we
establish the existence of a skew multiplicative $K$-basis for every
skew 2-nomial algebra, and we explore the existence of a (left,
right, or two-sided) monomial ordering for an (almost) skew 2-nomial
algebra. As distinct from commonly recognized algebras holding a
Gr\"obner basis theory (such as algebras of the solvable type [K-RW]
and some of their homomorphic images), a subclass of skew 2-nomial
algebras that have a left  Gr\"obner basis theory but may not
necessarily have a two-sided Gr\"obner basis theory, respectively a
subclass of skew 2-nomial algebras that have a right  Gr\"obner
basis theory but may not necessarily have a two-sided Gr\"obner
basis theory, are determined such that numerous quantum binomial
algebras (which provide binomial solutions to the Yang-baxter
equation [Laf], [G-I2]) are involved. }
\end{minipage}\end{center}

\def\QED{\hfill{$\Box$}} \def\NZ{\mathbb{N}}
\def \r{\rightarrow}

\def\normalbaselines{\baselineskip 24pt\lineskip 4pt\lineskiplimit 4pt}
\def\mapdown#1{\llap{$\vcenter {\hbox {$\scriptstyle #1$}}$}
                                \Bigg\downarrow}
\def\mapdownr#1{\Bigg\downarrow\rlap{$\vcenter{\hbox
                                    {$\scriptstyle #1$}}$}}
\def\mapright#1#2{\smash{\mathop{\longrightarrow}\limits^{#1}_{#2}}}
\def\mapleft#1#2{\smash{\mathop{\longleftarrow}\limits^{#1}_{#2}}}
\def\mapup#1{\Bigg\uparrow\rlap{$\vcenter {\hbox  {$\scriptstyle #1$}}$}}
\def\mapupl#1{\llap{$\vcenter {\hbox {$\scriptstyle #1$}}$}
                                      \Bigg\uparrow}
\def\v5{\vskip .5truecm}
\def\T#1{\widetilde #1}
\def\OV#1{\overline {#1}}
\def\hang{\hangindent\parindent}
\def\textindent#1{\indent\llap{#1\enspace}\ignorespaces}
\def\item{\par\hang\textindent}
\message{<Paul Taylor's commutative diagrams, 20 July 1990>}
\newdimen\DiagramCellHeight\DiagramCellHeight3em 
\newdimen\DiagramCellWidth\DiagramCellWidth3em 
\newdimen\MapBreadth\MapBreadth.04em 
\newdimen\MapShortFall\MapShortFall.4em 
\newdimen\PileSpacing\PileSpacing1em 
\def\labelstyle{\ifincommdiag\textstyle\else\scriptstyle\fi}
\let\objectstyle\displaystyle


\def\rTo{\HorizontalMap\empty-\empty-\rhvee}
\def\lTo{\HorizontalMap\lhvee-\empty-\empty}
\def\dTo{\VerticalMap\empty|\empty|\dhvee}
\def\uTo{\VerticalMap\uhvee|\empty|\empty}
\let\uFrom\uTo\let\lFrom\lTo

\def\rArr{\HorizontalMap\empty-\empty-\rhla}
\def\lArr{\HorizontalMap\lhla-\empty-\empty}
\def\dArr{\VerticalMap\empty|\empty|\dhla}
\def\uArr{\VerticalMap\uhla|\empty|\empty}

\def\rDotsto{\HorizontalMap\empty\hfdot\hfdot\hfdot\rhvee}
\def\lDotsto{\HorizontalMap\lhvee\hfdot\hfdot\hfdot\empty}
\def\dDotsto{\VerticalMap\empty\vfdot\vfdot\vfdot\dhvee}
\def\uDotsto{\VerticalMap\uhvee\vfdot\vfdot\vfdot\empty}
\let\uDotsfrom\uDotsto\let\lDotsfrom\lDotsto

\def\rDashto{\HorizontalMap\empty\hfdash\hfdash\hfdash\rhvee}
\def\lDashto{\HorizontalMap\lhvee\hfdash\hfdash\hfdash\empty}
\def\dDashto{\VerticalMap\empty\vfdash\vfdash\vfdash\dhvee}
\def\uDashto{\VerticalMap\uhvee\vfdash\vfdash\vfdash\empty}
\let\uDashfrom\uDashto\let\lDashfrom\lDashto

\def\rImplies{\HorizontalMap==\empty=\Rightarrow}
\def\lImplies{\HorizontalMap\Leftarrow=\empty==}
\def\dImplies{\VerticalMap\|\|\empty\|\Downarrow}
\def\uImplies{\VerticalMap\Uparrow\|\empty\|\|}
\let\uImpliedby\uImplies\let\lImpliedby\lImplies

\def\rMapsto{\HorizontalMap\rtbar-\empty-\rhvee}
\def\lMapsto{\HorizontalMap\lhvee-\empty-\ltbar}
\def\dMapsto{\VerticalMap\dtbar|\empty|\dhvee}
\def\uMapsto{\VerticalMap\uhvee|\empty|\utbar}
\let\uMapsfrom\uMapsto\let\lMapsfrom\lMapsto

\def\rIntoA{\HorizontalMap\rthooka-\empty-\rhvee}
\def\rIntoB{\HorizontalMap\rthookb-\empty-\rhvee}
\def\lIntoA{\HorizontalMap\lhvee-\empty-\lthooka}
\def\lIntoB{\HorizontalMap\lhvee-\empty-\lthookb}
\def\dIntoA{\VerticalMap\dthooka|\empty|\dhvee}
\def\dIntoB{\VerticalMap\dthookb|\empty|\dhvee}
\def\uIntoA{\VerticalMap\uhvee|\empty|\uthooka}
\def\uIntoB{\VerticalMap\uhvee|\empty|\uthookb}
\let\uInfromA\uIntoA\let\uInfromB\uIntoB\let\lInfromA\lIntoA\let\lInfromB
\lIntoB\let\rInto\rIntoA\let\lInto\lIntoA\let\dInto\dIntoB\let\uInto\uIntoA

\def\rEmbed{\HorizontalMap\gt-\empty-\rhvee}
\def\lEmbed{\HorizontalMap\lhvee-\empty-\lt}
\def\dEmbed{\VerticalMap\vee|\empty|\dhvee}
\def\uEmbed{\VerticalMap\uhvee|\empty|\wedge}

\def\rProject{\HorizontalMap\empty-\empty-\triangleright}
\def\lProject{\HorizontalMap\triangleleft-\empty-\empty}
\def\uProject{\VerticalMap\triangleup|\empty|\empty}
\def\dProject{\VerticalMap\empty|\empty|\littletriangledown}

\def\rOnto{\HorizontalMap\empty-\empty-\twoheadrightarrow}
\def\lOnto{\HorizontalMap\twoheadleftarrow-\empty-\empty}
\def\dOnto{\VerticalMap\empty|\empty|\twoheaddownarrow}
\def\uOnto{\VerticalMap\twoheaduparrow|\empty|\empty}
\let\lOnfrom\lOnto\let\uOnfrom\uOnto

\def\hEq{\HorizontalMap==\empty==}
\def\vEq{\VerticalMap\|\|\empty\|\|}
\let\rEq\hEq\let\lEq\hEq\let\uEq\vEq\let\dEq\vEq

\def\hLine{\HorizontalMap\empty-\empty-\empty}
\def\vLine{\VerticalMap\empty|\empty|\empty}
\let\rLine\hLine\let\lLine\hLine\let\uLine\vLine\let\dLine\vLine

\def\hDots{\HorizontalMap\empty\hfdot\hfdot\hfdot\empty}
\def\vDots{\VerticalMap\empty\vfdot\vfdot\vfdot\empty}
\let\rDots\hDots\let\lDots\hDots\let\uDots\vDots\let\dDots\vDots

\def\hDashes{\HorizontalMap\empty\hfdash\hfdash\hfdash\empty}
\def\vDashes{\VerticalMap\empty\vfdash\vfdash\vfdash\empty}
\let\rDashes\hDashes\let\lDashes\hDashes\let\uDashes\vDashes\let\dDashes
\vDashes

\def\rPto{\HorizontalMap\empty-\empty-\rightharpoonup}
\def\lPto{\HorizontalMap\leftharpoonup-\empty-\empty}
\def\uPto{\VerticalMap\upharpoonright|\empty|\empty}
\def\dPto{\VerticalMap\empty|\empty|\downharpoonright}
\let\lPfrom\lPto\let\uPfrom\uPto

\def\NW{\NorthWest\DiagonalMap{\lah111}{\laf100}{}{\laf100}{}(2,2)}
\def\NE{\NorthEast\DiagonalMap{\lah22}{\laf0}{}{\laf0}{}(2,2)}
\def\SW{\SouthWest\DiagonalMap{}{\laf0}{}{\laf0}{\lah11}(2,2)}
\def\SE{\SouthEast\DiagonalMap{}{\laf100}{}{\laf100}{\lah122}(2,2)}

\def\nNW{\NorthWest\DiagonalMap{\lah135}{\laf112}{}{\laf112}{}(2,3)}
\def\nNE{\NorthEast\DiagonalMap{\lah36}{\laf12}{}{\laf12}{}(2,3)}
\def\sSW{\SouthWest\DiagonalMap{}{\laf12}{}{\laf12}{\lah35}(2,3)}
\def\sSE{\SouthEast\DiagonalMap{}{\laf112}{}{\laf112}{\lah136}(2,3)}

\def\wNW{\NorthWest\DiagonalMap{\lah153}{\laf121}{}{\laf121}{}(3,2)}
\def\eNE{\NorthEast\DiagonalMap{\lah63}{\laf21}{}{\laf21}{}(3,2)}
\def\wSW{\SouthWest\DiagonalMap{}{\laf21}{}{\laf21}{\lah53}(3,2)}
\def\eSE{\SouthEast\DiagonalMap{}{\laf121}{}{\laf121}{\lah163}(3,2)}

\def\NNW{\NorthWest\DiagonalMap{\lah113}{\laf101}{}{\laf101}{}(2,4)}
\def\NNE{\NorthEast\DiagonalMap{\lah25}{\laf01}{}{\laf01}{}(2,4)}
\def\SSW{\SouthWest\DiagonalMap{}{\laf01}{}{\laf01}{\lah13}(2,4)}
\def\SSE{\SouthEast\DiagonalMap{}{\laf101}{}{\laf101}{\lah125}(2,4)}

\def\WNW{\NorthWest\DiagonalMap{\lah131}{\laf110}{}{\laf110}{}(4,2)}
\def\ENE{\NorthEast\DiagonalMap{\lah52}{\laf10}{}{\laf10}{}(4,2)}
\def\WSW{\SouthWest\DiagonalMap{}{\laf10}{}{\laf10}{\lah31}(4,2)}
\def\ESE{\SouthEast\DiagonalMap{}{\laf110}{}{\laf110}{\lah152}(4,2)}

\def\NNNW{\NorthWest\DiagonalMap{\lah115}{\laf102}{}{\laf102}{}(2,6)}
\def\NNNE{\NorthEast\DiagonalMap{\lah16}{\laf02}{}{\laf02}{}(2,6)}
\def\SSSW{\SouthWest\DiagonalMap{}{\laf02}{}{\laf02}{\lah15}(2,6)}
\def\SSSE{\SouthEast\DiagonalMap{}{\laf102}{}{\laf102}{\lah116}(2,6)}

\def\WWNW{\NorthWest\DiagonalMap{\lah151}{\laf120}{}{\laf120}{}(6,2)}
\def\EENE{\NorthEast\DiagonalMap{\lah61}{\laf20}{}{\laf20}{}(6,2)}
\def\WWSW{\SouthWest\DiagonalMap{}{\laf20}{}{\laf20}{\lah51}(6,2)}
\def\EESE{\SouthEast\DiagonalMap{}{\laf120}{}{\laf120}{\lah161}(6,2)}


\font\tenln=line10

\mathchardef\lt="313C \mathchardef\gt="313E

\def\rhvee{\mkern-10mu\gt}
\def\lhvee{\lt\mkern-10mu}
\def\dhvee{\vbox\tozpt{\vss\hbox{$\vee$}\kern0pt}}
\def\uhvee{\vbox\tozpt{\hbox{$\wedge$}\vss}}
\def\rhcvee{\mkern-10mu\succ}
\def\lhcvee{\prec\mkern-10mu}
\def\dhcvee{\vbox\tozpt{\vss\hbox{$\curlyvee$}\kern0pt}}
\def\uhcvee{\vbox\tozpt{\hbox{$\curlywedge$}\vss}}
\def\rhvvee{\mkern-10mu\gg}
\def\lhvvee{\ll\mkern-10mu}
\def\dhvvee{\vbox\tozpt{\vss\hbox{$\vee$}\kern-.6ex\hbox{$\vee$}\kern0pt}}
\def\uhvvee{\vbox\tozpt{\hbox{$\wedge$}\kern-.6ex\hbox{$\wedge$}\vss}}
\def\twoheaddownarrow{\rlap{$\downarrow$}\raise-.5ex\hbox{$\downarrow$}}
\def\twoheaduparrow{\rlap{$\uparrow$}\raise.5ex\hbox{$\uparrow$}}
\def\triangleup{{\scriptscriptstyle\bigtriangleup}}
\def\littletriangledown{{\scriptscriptstyle\triangledown}}
\def\rhla{\vbox\tozpt{\vss\hbox\tozpt{\hss\tenln\char'55}\kern\axisheight}}
\def\lhla{\vbox\tozpt{\vss\hbox\tozpt{\tenln\char'33\hss}\kern\axisheight}}
\def\dhla{\vbox\tozpt{\vss\hbox\tozpt{\tenln\char'77\hss}}}
\def\uhla{\vbox\tozpt{\hbox\tozpt{\tenln\char'66\hss}\vss}}
\def\htdot{\mkern3.15mu\cdot\mkern3.15mu}
\def\vtdot{\vbox to 1.46ex{\vss\hbox{$\cdot$}}}
\def\utbar{\vrule height 0.093ex depth0pt width 0.4em} \let\dtbar\utbar
\def\rtbar{\mkern1.5mu\vrule height 1.1ex depth.06ex width .04em\mkern1.5mu}%
\let\ltbar\rtbar
\def\rthooka{\raise.603ex\hbox{$\scriptscriptstyle\subset$}}
\def\lthooka{\raise.603ex\hbox{$\scriptscriptstyle\supset$}}
\def\rthookb{\raise-.022ex\hbox{$\scriptscriptstyle\subset$}}
\def\lthookb{\raise-.022ex\hbox{$\scriptscriptstyle\supset$}}
\def\dthookb{\hbox{$\scriptscriptstyle\cap$}\mkern5.5mu}
\def\uthookb{\hbox{$\scriptscriptstyle\cup$}\mkern4.5mu}
\def\dthooka{\mkern6mu\hbox{$\scriptscriptstyle\cap$}}
\def\uthooka{\mkern6mu\hbox{$\scriptscriptstyle\cup$}}
\def\hfdot{\mkern3.15mu\cdot\mkern3.15mu}
\def\vfdot{\vbox to 1.46ex{\vss\hbox{$\cdot$}}}
\def\vfdashstrut{\vrule width0pt height1.3ex depth0.7ex}
\def\vfthedash{\vrule width\MapBreadth height0.6ex depth 0pt}
\def\hfthedash{\vrule\horizhtdp width 0.26em}
\def\hfdash{\mkern5.5mu\hfthedash\mkern5.5mu}
\def\vfdash{\vfdashstrut\vfthedash}

\def\nwhTO{\nwarrow\mkern-1mu}
\def\nehTO{\mkern-.1mu\nearrow}
\def\sehTO{\searrow\mkern-.02mu}
\def\swhTO{\mkern-.8mu\swarrow}

\def\SEpbk{\rlap{\smash{\kern0.1em \vrule depth 2.67ex height -2.55ex width 0%
.9em \vrule height -0.46ex depth 2.67ex width .05em }}}
\def\SWpbk{\llap{\smash{\vrule height -0.46ex depth 2.67ex width .05em \vrule
depth 2.67ex height -2.55ex width .9em \kern0.1em }}}
\def\NEpbk{\rlap{\smash{\kern0.1em \vrule depth -3.48ex height 3.67ex width 0%
.95em \vrule height 3.67ex depth -1.39ex width .05em }}}
\def\NWpbk{\llap{\smash{\vrule height 3.6ex depth -1.39ex width .05em \vrule
depth -3.48ex height 3.67ex width .95em \kern0.1em }}}


\newcount\cdna\newcount\cdnb\newcount\cdnc\newcount\cdnd\cdna=\catcode`\@%
\catcode`\@=11 \let\then\relax\def\loopa#1\repeat{\def\bodya{#1}\iteratea}%
\def\iteratea{\bodya\let\next\iteratea\else\let\next\relax\fi\next}\def\loopb
#1\repeat{\def\bodyb{#1}\iterateb}\def\iterateb{\bodyb\let\next\iterateb\else
\let\next\relax\fi\next} \def\mapctxterr{\message{commutative diagram: map
context error}}\def\mapclasherr{\message{commutative diagram: clashing maps}}%
\def\ObsDim#1{\expandafter\message{! diagrams Warning: Dimension \string#1 is
obsolete
(ignored)}\global\let#1\ObsDimq\ObsDimq}\def\ObsDimq{\dimen@=}\def
\HorizontalMapLength{\ObsDim\HorizontalMapLength}\def\VerticalMapHeight{%
\ObsDim\VerticalMapHeight}\def\VerticalMapDepth{\ObsDim\VerticalMapDepth}\def
\VerticalMapExtraHeight{\ObsDim\VerticalMapExtraHeight}\def
\VerticalMapExtraDepth{\ObsDim\VerticalMapExtraDepth}\def\ObsCount#1{%
\expandafter\message{! diagrams Warning: Count \string#1 is obsolete (ignored%
)}\global\let#1\ObsCountq\ObsCountq}\def\ObsCountq{\count@=}\def
\DiagonalLineSegments{\ObsCount\DiagonalLineSegments}\def\tozpt{to\z@}\def
\sethorizhtdp{\dimen8=\axisheight\dimen9=\MapBreadth\advance\dimen8.5\dimen9%
\advance\dimen9-\dimen8}\def\horizhtdp{height\dimen8 depth\dimen9
}\def \axisheight{\fontdimen22\the\textfont2 }\countdef\boxc@unt=14

\def\bombparameters{\hsize\z@\rightskip\z@ plus1fil minus\maxdimen
\parfillskip\z@\linepenalty9000 \looseness0 \hfuzz\maxdimen\hbadness10000
\clubpenalty0 \widowpenalty0 \displaywidowpenalty0
\interlinepenalty0 \predisplaypenalty0 \postdisplaypenalty0
\interdisplaylinepenalty0
\interfootnotelinepenalty0 \floatingpenalty0 \brokenpenalty0 \everypar{}%
\leftskip\z@\parskip\z@\parindent\z@\pretolerance10000
\tolerance10000 \hyphenpenalty10000 \exhyphenpenalty10000
\binoppenalty10000 \relpenalty10000 \adjdemerits0
\doublehyphendemerits0 \finalhyphendemerits0 \prevdepth\z@}\def
\startbombverticallist{\hbox{}\penalty1\nointerlineskip}

\def\pushh#1\to#2{\setbox#2=\hbox{\box#1\unhbox#2}}\def\pusht#1\to#2{\setbox#%
2=\hbox{\unhbox#2\box#1}}

\newif\ifallowhorizmap\allowhorizmaptrue\newif\ifallowvertmap
\allowvertmapfalse\newif\ifincommdiag\incommdiagfalse

\def\diagram{\hbox\bgroup$\vcenter\bgroup\startbombverticallist
\incommdiagtrue\baselineskip\DiagramCellHeight\lineskip\z@\lineskiplimit\z@
\mathsurround\z@\tabskip\z@\let\\\diagcr\allowhorizmaptrue\allowvertmaptrue
\halign\bgroup\lcdtempl##\rcdtempl&&\lcdtempl##\rcdtempl\cr}\def\enddiagram{%
\crcr\egroup\reformatmatrix\egroup$\egroup}\def\commdiag#1{{\diagram#1%
\enddiagram}}

\def\lcdtempl{\futurelet\thefirsttoken\dolcdtempl}\newif\ifemptycell\def
\dolcdtempl{\ifx\thefirsttoken\rcdtempl\then\hskip1sp plus 1fil
\emptycelltrue
\else\hfil$\emptycellfalse\objectstyle\fi}\def\rcdtempl{\ifemptycell\else$%
\hfil\fi}\def\diagcr{\cr} \def\across#1{\span\omit\mscount=#1
\loop\ifnum \mscount>2
\spAn\repeat\ignorespaces}\def\spAn{\relax\span\omit\advance
\mscount by -1}

\def\CellSize{\afterassignment\cdhttowd\DiagramCellHeight}\def\cdhttowd{%
\DiagramCellWidth\DiagramCellHeight}\def\MapsAbut{\MapShortFall\z@}

\newcount\cdvdl\newcount\cdvdr\newcount\cdvd\newcount\cdbfb\newcount\cdbfr
\newcount\cdbfl\newcount\cdvdr\newcount\cdvdl\newcount\cdvd

\def\reformatmatrix{\bombparameters\cdvdl=\insc@unt\cdvdr=\cdvdl\cdbfb=%
\boxc@unt\advance\cdbfb1
\cdbfr=\cdbfb\setbox1=\vbox{}\dimen2=\z@\loop\setbox
0=\lastbox\ifhbox0 \dimen1=\lastskip\unskip\dimen5=\ht0
\advance\dimen5 \dimen
1 \dimen4=\dp0 \penalty1 \reformatrow\unpenalty\ht4=\dimen5 \dp4=\dimen4 \ht3%
\z@\dp3\z@\setbox1=\vbox{\box4 \nointerlineskip\box3 \nointerlineskip\unvbox1%
}\dimen2=\dimen1 \repeat\unvbox1}

\newif\ifcontinuerow

\def\reformatrow{\cdbfl=\cdbfr\noindent\unhbox0 \loopa\unskip\setbox\cdbfl=%
\lastbox\ifhbox\cdbfl\advance\cdbfl1\repeat\par\unskip\dimen6=2%
\DiagramCellWidth\dimen7=-\DiagramCellWidth\setbox3=\hbox{}\setbox4=\hbox{}%
\setbox7=\box\voidb@x\cdvd=\cdvdl\continuerowtrue\loopa\advance\cdvd-1
\adjustcells\ifcontinuerow\advance\dimen6\wd\cdbfl\cdda=.5\dimen6
\ifdim\cdda
<\DiagramCellWidth\then\dimen6\DiagramCellWidth\advance\dimen6-\cdda
\nopendvert\cdda\DiagramCellWidth\fi\advance\dimen7\cdda\dimen6=\wd\cdbfl
\reformatcell\advance\cdbfl-1 \repeat\advance\dimen7.5\dimen6
\outHarrow} \def
\adjustcells{\ifnum\cdbfr>\cdbfl\then\ifnum\cdvdr>\cdvd\then\continuerowfalse
\else\setbox\cdbfl=\hbox
to\wd\cdvd{\lcdtempl\VonH{}\rcdtempl}\fi\else\ifnum
\cdvdr>\cdvd\then\advance\cdvdr-1
\setbox\cdvd=\vbox{}\wd\cdvd=\wd\cdbfl\dp \cdvd=\dp1 \fi\fi}

\def\reformatcell{\sethorizhtdp\noindent\unhbox\cdbfl\skip0=\lastskip\unskip
\par\ifcase\prevgraf\reformatempty\or\reformatobject\else\reformatcomplex\fi
\unskip}\def\reformatobject{\setbox6=\lastbox\unskip\vadjdon6\outVarrow
\setbox6=\hbox{\unhbox6}\advance\dimen7-.5\wd6
\outHarrow\dimen7=-.5\wd6 \pusht6\to4}\newcount\globnum

\def\reformatcomplex{\setbox6=\lastbox\unskip\setbox9=\lastbox\unskip\setbox9%
=\hbox{\unhbox9
\skip0=\lastskip\unskip\global\globnum\lastpenalty\hskip\skip 0
}\advance\globnum9999
\ifcase\globnum\reformathoriz\or\reformatpile\or
\reformatHonV\or\reformatVonH\or\reformatvert\or\reformatHmeetV\fi}

\def\reformatempty{\vpassdon\ifdim\skip0>\z@\then\hpassdon\else\ifvoid2 \then
\else\advance\dimen7-.5\dimen0
\cdda=\wd2\advance\cdda.5\dimen0\wd2=\cdda\fi
\fi}\def\VonH{\doVonH6}\def\HonV{\doVonH7}\def\HmeetV{\MapBreadth-2%
\MapShortFall\doVonH4}\def\doVonH#1{\cdna-999#1\futurelet\thenexttoken
\dooVonH}\def\dooVonH{\let\next\relax\sethorizhtdp\ifallowhorizmap
\ifallowvertmap\then\ifx\thenexttoken[\then\let\next\VonHstrut\else
\sethorizhtdp\dimen0\MapBreadth\let\next\VonHnostrut\fi\else\mapctxterr\fi
\else\mapctxterr\fi\next}\def\VonHstrut[#1]{\setbox0=\hbox{$#1$}\dimen0\wd0%
\dimen8\ht0\dimen9\dp0 \VonHnostrut}\def\VonHnostrut{\setbox0=\hbox{}\ht0=%
\dimen8\dp0=\dimen9\wd0=.5\dimen0 \copy0\penalty\cdna\box0
\allowhorizmapfalse
\allowvertmapfalse}\def\reformatHonV{\hpassdon\doreformatHonV}\def
\reformatHmeetV{\dimen@=\wd9 \advance\dimen7-\wd9 \outHarrow\setbox6=\hbox{%
\unhbox6}\dimen7-\wd6 \advance\dimen@\wd6 \setbox6=\hbox to\dimen@{\hss}%
\pusht6\to4\doreformatHonV}\def\doreformatHonV{\setbox9=\hbox{\unhbox9
\unskip
\unpenalty\global\setbox\globbox=\lastbox}\vadjdon\globbox\outVarrow}\def
\reformatVonH{\vpassdon\advance\dimen7-\wd9 \outHarrow\setbox6=\hbox{\unhbox6%
}\dimen7=-\wd6 \setbox6=\hbox{\kern\wd9 \kern\wd6}\pusht6\to4}\def\hpassdon{}%
\def\vpassdon{\dimen@=\dp\cdvd\advance\dimen@\dimen4 \advance\dimen@\dimen5
\dp\cdvd=\dimen@\nopendvert}\def\vadjdon#1{\dimen8=\ht#1
\dimen9=\dp#1 }

\def\HorizontalMap#1#2#3#4#5{\sethorizhtdp\setbox1=\makeharrowpart{#1}\def
\arrowfillera{#2}\def\arrowfillerb{#4}\setbox5=\makeharrowpart{#5}\ifx
\arrowfillera\justhorizline\then\def\arra{\hrule\horizhtdp}\def\kea{\kern-0.%
01em}\let\arrstruthtdp\horizhtdp\else\def\kea{\kern-0.15em}\setbox2=\hbox{%
\kea${\arrowfillera}$\kea}\def\arra{\copy2}\def\arrstruthtdp{height\ht2 depth%
\dp2
}\fi\ifx\arrowfillerb\justhorizline\then\def\arrb{\hrule\horizhtdp}\def
\keb{kern-0.01em}\ifx\arrowfillera\empty\then\let\arrstruthtdp\horizhtdp\fi
\else\def\keb{\kern-0.15em}\setbox4=\hbox{\keb${\arrowfillerb}$\keb}\def\arrb
{\copy4}\ifx\arrowfilera\empty\then\def\arrstruthtdp{height\ht4 depth\dp4 }%
\fi\fi\setbox3=\makeharrowpart{{#3}\vrule width\z@\arrstruthtdp}%
\ifallowhorizmap\then\let\execmap\execHorizontalMap\else\let\execmap
\mapctxterr\fi\allowhorizmapfalse\gettwoargs}\def\makeharrowpart#1{\hbox{%
\mathsurround\z@\edef\next{#1}\ifx\next\empty\else$\mkern-1.5mu{\next}\mkern-%
1.5mu$\fi}}\def\justhorizline{-}

\def\execHorizontalMap{\dimen0=\wd6 \ifdim\dimen0<\wd7\then\dimen0=\wd7\fi
\dimen3=\wd3 \ifdim\dimen0<2em\then\dimen0=2em\fi\skip2=.5\dimen0
\ifincommdiag plus 1fill\fi minus\z@\advance\skip2-.5\dimen3
\skip4=\skip2 \advance\skip2-\wd1 \advance\skip4-\wd5
\kern\MapShortFall\box1 \xleaders
\arra\hskip\skip2 \vbox{\lineskiplimit\maxdimen\lineskip.5ex \ifhbox6 \hbox to%
\dimen3 {\hss\box6\hss}\fi\vtop{\box3 \ifhbox7 \hbox to\dimen3
{\hss\box7\hss
}\fi}}\ifincommdiag\kern-.5\dimen3\penalty-9999\null\kern.5\dimen3\fi
\xleaders\arrb\hskip\skip4 \box5 \kern\MapShortFall}

\def\reformathoriz{\vadjdon6\outVarrow\ifvoid7\else\mapclasherr\fi\setbox2=%
\box9 \wd2=\dimen7 \dimen7=\z@\setbox7=\box6 }

\def\resetharrowpart#1#2{\ifvoid#1\then\ifdim#2=\z@\else\setbox4=\hbox{%
\unhbox4\kern#2}\fi\else\ifhbox#1\then\setbox#1=\hbox
to#2{\unhbox#1}\else
\widenpile#1\fi\pusht#1\to4\fi}\def\outHarrow{\resetharrowpart2{\wd2}\pusht2%
\to4\resetharrowpart7{\dimen7}\pusht7\to4\dimen7=\z@}

\def\pile#1{{\incommdiagtrue\let\pile\innerpile\allowvertmapfalse
\allowhorizmaptrue\baselineskip.5\PileSpacing\lineskip\z@\lineskiplimit\z@
\mathsurround\z@\tabskip\z@\let\\\pilecr\vcenter{\halign{\hfil$##$\hfil\cr#1
\crcr}}}\ifincommdiag\then\ifallowhorizmap\then\penalty-9998
\allowvertmapfalse\allowhorizmapfalse\else\mapctxterr\fi\fi}\def\pilecr{\cr}%
\def\innerpile#1{\noalign{\halign{\hfil$##$\hfil\cr#1 \crcr}}}

\def\reformatpile{\vadjdon9\outVarrow\ifvoid7\else\mapclasherr\fi\penalty1
\setbox9=\hbox{\unhbox9 \unskip\unpenalty\setbox9=\lastbox\unhbox9
\global
\setbox\globbox=\lastbox}\unvbox\globbox\setbox9=\vbox{}\setbox7=\vbox{}%
\loopb\setbox6=\lastbox\ifhbox6
\skip3=\lastskip\unskip\splitpilerow\repeat
\unpenalty\setbox9=\hbox{$\vcenter{\unvbox9}$}\setbox2=\box9
\dimen7=\z@}\def
\pilestrut{\vrule height\dimen0 depth\dimen3 width\z@}\def\splitpilerow{%
\dimen0=\ht6 \dimen3=\dp6
\noindent\unhbox6\unskip\setbox6=\lastbox\unskip
\unhbox6\par\setbox6=\lastbox\unskip\ifcase\prevgraf\or\setbox6=\hbox\tozpt{%
\hss\unhbox6\hss}\ht6=\dimen0 \dp6=\dimen3
\setbox9=\vbox{\vskip\skip3 \hbox
to\dimen7{\hfil\box6}\nointerlineskip\unvbox9}\setbox7=\vbox{\vskip\skip3
\hbox{\pilestrut\hfil}\nointerlineskip\unvbox7}\or\setbox7=\vbox{\vskip\skip3
\hbox{\pilestrut\unhbox6}\nointerlineskip\unvbox7}\setbox6=\lastbox\unskip
\setbox9=\vbox{\vskip\skip3 \hbox to\dimen7{\pilestrut\unhbox6}%
\nointerlineskip\unvbox9}\fi\unskip}

\def\widenpile#1{\setbox#1=\hbox{$\vcenter{\unvbox#1 \setbox8=\vbox{}\loopb
\setbox9=\lastbox\ifhbox9
\skip3=\lastskip\unskip\setbox8=\vbox{\vskip\skip3 \hbox
to\dimen7{\unhbox9}\nointerlineskip\unvbox8}\repeat\unvbox8 }$}}

\def\justverticalline{|}\def\makevarrowpart#1{\hbox to\MapBreadth{\hss$\kern
\MapBreadth{#1}$\hss}}\def\VerticalMap#1#2#3#4#5{\setbox1=\makevarrowpart{#1}%
\def\arrowfillera{#2}\setbox3=\makevarrowpart{#3}\def\arrowfillerb{#4}\setbox
5=\makevarrowpart{#5}\ifx\arrowfillera\justverticalline\then\def\arra{\vrule
width\MapBreadth}\def\kea{\kern-0.05ex}\else\def\kea{\kern-0.35ex}\setbox2=%
\vbox{\kea\makevarrowpart\arrowfillera\kea}\def\arra{\copy2}\fi\ifx
\arrowfillerb\justverticalline\then\def\arrb{\vrule
width\MapBreadth}\def\keb
{\kern-0.05ex}\else\def\keb{\kern-0.35ex}\setbox4=\vbox{\keb\makevarrowpart
\arrowfillerb\keb}\def\arrb{\copy4}\fi\ifallowvertmap\then\let\execmap
\execVerticalMap\else\let\execmap\mapctxterr\fi\allowhorizmapfalse\gettwoargs
}

\def\execVerticalMap{\setbox3=\makevarrowpart{\box3}\setbox0=\hbox{}\ht0=\ht3
\dp0\z@\ht3\z@\box6 \setbox8=\vtop spread2ex{\offinterlineskip\box3
\xleaders
\arrb\vfill\box5 \kern\MapShortFall}\dp8=\z@\box8 \kern-\MapBreadth\setbox8=%
\vbox spread2ex{\offinterlineskip\kern\MapShortFall\box1
\xleaders\arra\vfill \box0}\ht8=\z@\box8
\ifincommdiag\then\kern-.5\MapBreadth\penalty-9995 \null
\kern.5\MapBreadth\fi\box7\hfil}

\newcount\colno\newdimen\cdda\newbox\globbox\def\reformatvert{\setbox6=\hbox{%
\unhbox6}\cdda=\wd6 \dimen3=\dp\cdvd\advance\dimen3\dimen4
\setbox\cdvd=\hbox {}\colno=\prevgraf\advance\colno-2
\loopb\setbox9=\hbox{\unhbox9 \unskip
\unpenalty\dimen7=\lastkern\unkern\global\setbox\globbox=\lastbox\advance
\dimen7\wd\globbox\advance\dimen7\lastkern\unkern\setbox9=\lastbox\vtop to%
\dimen3{\unvbox9}\kern\dimen7 }\ifnum\colno>0
\ifdim\wd9<\PileSpacing\then \setbox9=\hbox
to\PileSpacing{\unhbox9}\fi\dimen0=\wd9 \advance\dimen0-\wd
\globbox\setbox\cdvd=\hbox{\kern\dimen0 \box\globbox\unhbox\cdvd}\pushh9\to6%
\advance\colno-1
\setbox9=\lastbox\unskip\repeat\advance\dimen7-.5\wd6
\advance\dimen7.5\cdda\advance\dimen7-\wd9 \outHarrow\dimen7=-.5\wd6
\advance
\dimen7-.5\cdda\pusht9\to4\pusht6\to4\nopendvert\dimen@=\dimen6\advance
\dimen@-\wd\cdvd\advance\dimen@-\wd\globbox\divide\dimen@2
\setbox\cdvd=\hbox
{\kern\dimen@\box\globbox\unhbox\cdvd\kern\dimen@}\dimen8=\dp\cdvd\advance
\dimen8\dimen5 \dp\cdvd=\dimen8 \ht\cdvd=\z@}

\def\outVarrow{\ifhbox\cdvd\then\deepenbox\cdvd\pusht\cdvd\to3\else
\nopendvert\fi\dimen3=\dimen5 \advance\dimen3-\dimen8 \setbox\cdvd=\vbox{%
\vfil}\dp\cdvd=\dimen3} \def\nopendvert{\setbox3=\hbox{\unhbox3\kern\dimen6}}%
\def\deepenbox\cdvd{\setbox\cdvd=\hbox{\dimen3=\dimen4 \advance\dimen3-\dimen
9 \setbox6=\hbox{}\ht6=\dimen3 \dp6=-\dimen3 \dimen0=\dp\cdvd\advance\dimen0%
\dimen3 \unhbox\cdvd\dimen3=\lastkern\unkern\setbox8=\hbox{\kern\dimen3}%
\loopb\setbox9=\lastbox\ifvbox9 \setbox9=\vtop to\dimen0{\copy6
\nointerlineskip\unvbox9 }\dimen3=\lastkern\unkern\setbox8=\hbox{\kern\dimen3%
\box9\unhbox8}\repeat\unhbox8 }}

\newif\ifPositiveGradient\PositiveGradienttrue\newif\ifClimbing\Climbingtrue
\newcount\DiagonalChoice\DiagonalChoice1 \newcount\lineno\newcount\rowno
\newcount\charno\def\laf{\afterassignment\xlaf\charno='}\def\xlaf{\hbox{%
\tenln\char\charno}}\def\lah{\afterassignment\xlah\charno='}\def\xlah{\hbox{%
\tenln\char\charno}}\def\makedarrowpart#1{\hbox{\mathsurround\z@${#1}$}}\def
\lad{\afterassignment\xlad\charno='}\def\xlad{\setbox2=\xlaf\setbox0=\hbox to%
.5\wd2{$\hss\ldot\hss$}\ht0=.25\ht2 \dp0=\ht0 \hbox{\mv-\ht0\copy0 \mv\ht0%
\box0}}

\def\DiagonalMap#1#2#3#4#5{\ifPositiveGradient\then\let\mv\raise\else\let\mv
\lower\fi\setbox2=\makedarrowpart{#2}\setbox1=\makedarrowpart{#1}\setbox4=%
\makedarrowpart{#4}\setbox5=\makedarrowpart{#5}\setbox3=\makedarrowpart{#3}%
\let\execmap\execDiagonalLine\gettwoargs}

\def\makeline#1(#2,#3;#4){\hbox{\dimen1=#2\relax\dimen2=#3\relax\dimen5=#4%
\relax\vrule height\dimen5 depth\z@ width\z@\setbox8=\hbox to\dimen1{\tenln#1%
\hss}\cdna=\dimen5 \divide\cdna\dimen2 \ifnum\cdna=0 \then\box8 \else\dimen4=%
\dimen5 \advance\dimen4-\dimen2 \divide\dimen4\cdna\dimen3=\dimen1 \cdnb=%
\dimen2 \divide\cdnb1000 \divide\dimen3\cdnb\cdnb=\dimen4
\divide\cdnb1000 \multiply\dimen3\cdnb\dimen6\dimen1
\advance\dimen6-\dimen3 \cdnb=0
\ifPositiveGradient\then\dimen7\z@\else\dimen7\cdna\dimen4
\multiply\dimen4-1 \fi\loop\raise\dimen7\copy8
\ifnum\cdnb<\cdna\hskip-\dimen6 \advance\cdnb1
\advance\dimen7\dimen4
\repeat\fi}}\newdimen\objectheight\objectheight1.5ex

\def\execDiagonalLine{\setbox0=\hbox\tozpt{\cdna=\xcoord\cdnb=\ycoord\dimen8=%
\wd2 \dimen9=\ht2 \dimen0=\cdnb\DiagramCellHeight\advance\dimen0-2%
\MapShortFall\advance\dimen0-\objectheight\setbox2=\makeline\box2(\dimen8,%
\dimen9;.5\dimen0)\setbox4=\makeline\box4(\dimen8,\dimen9;.5\dimen0)\dimen0=2%
\wd2 \advance\dimen0-\cdna\DiagramCellWidth\advance\dimen0
2\DiagramCellWidth
\dimen2\DiagramCellHeight\advance\dimen2-\MapShortFall\dimen1\dimen2
\advance \dimen1-\ht1 \advance\dimen2-\ht2 \dimen6=\dimen2
\advance\dimen6.25\dimen8 \dimen3\dimen2 \advance\dimen3-\ht3
\dimen4=\dimen2 \dimen7=\dimen2 \advance \dimen4-\ht4
\advance\dimen7-\ht7 \advance\dimen7-.25\dimen8
\ifPositiveGradient\then\hss\raise\dimen4\hbox{\rlap{\box5}\box4}\llap{\raise
\dimen6\box6\kern.25\dimen9}\else\kern-.5\dimen0 \rlap{\raise\dimen1\box1}%
\raise\dimen2\box2 \llap{\raise\dimen7\box7\kern.25\dimen9}\fi\raise\dimen3%
\hbox\tozpt{\hss\box3\hss}\ifPositiveGradient\then\rlap{\kern.25\dimen9\raise
\dimen7\box7}\raise\dimen2\box2\llap{\raise\dimen1\box1}\kern-.5\dimen0
\else
\rlap{\kern.25\dimen9\raise\dimen6\box6}\raise\dimen4\hbox{\box4\llap{\box5}}%
\hss\fi}\ht0\z@\dp0\z@\box0}

\def\NorthWest{\PositiveGradientfalse\Climbingtrue\DiagonalChoice0 }\def
\NorthEast{\PositiveGradienttrue\Climbingtrue\DiagonalChoice1
}\def\SouthWest
{\PositiveGradienttrue\Climbingfalse\DiagonalChoice3 }\def\SouthEast{%
\PositiveGradientfalse\Climbingfalse\DiagonalChoice2 }

\newif\ifmoremapargs\def\gettwoargs{\setbox7=\box\voidb@x\setbox6=\box
\voidb@x\moremapargstrue\def\whichlabel{6}\def\xcoord{2}\def\ycoord{2}\def
\contgetarg{\def\whichlabel{7}\ifmoremapargs\then\let\next\getanarg\let
\contgetarg\execmap\else\let\next\execmap\fi\next}\getanarg}\def\getanarg{%
\futurelet\thenexttoken\switcharg}\def\getlabel#1#2#3{\setbox#1=\hbox{$%
\labelstyle\>{#3}\>$}\dimen0=\ht#1\advance\dimen0 .4ex\ht#1=\dimen0 \dimen0=%
\dp#1\advance\dimen0 .4ex\dp#1=\dimen0 \contgetarg}\def\eatspacerepeat{%
\afterassignment\getanarg\let\junk=
}\def\catcase#1:{{\ifcat\noexpand
\thenexttoken#1\then\global\let\xcase\docase\fi}\xcase}\def\tokcase#1:{{\ifx
\thenexttoken#1\then\global\let\xcase\docase\fi}\xcase}\def\default:{\docase}%
\def\docase#1\esac#2\esacs{#1}\def\skipcase#1\esac{}\def\getcoordsrepeat(#1,#%
2){\def\xcoord{#1}\def\ycoord{#2}\getanarg}\let\esacs\relax\def\switcharg{%
\global\let\xcase\skipcase\catcase{&}:\moremapargsfalse\contgetarg\esac
\catcase\bgroup:\getlabel\whichlabel-\esac\catcase^:\getlabel6\esac\catcase_:%
\getlabel7\esac\tokcase{~}:\getlabel3\esac\tokcase(:\getcoordsrepeat\esac
\catcase{
}:\eatspacerepeat\esac\default:\moremapargsfalse\contgetarg\esac
\esacs}

\catcode`\@=\cdna

\def\LH{{\bf LH}}\def\LM{{\bf LM}}\def\LT{{\bf
LT}}\def\KS{K\langle X\rangle}
\def\B{{\cal B}} \def\LC{{\bf LC}} \def\G{{\cal G}} \def\FRAC#1#2{\displaystyle{\frac{#1}{#2}}}
\def\SUM^#1_#2{\displaystyle{\sum^{#1}_{#2}}} \def\O{{\cal O}}  \def\J{{\bf J}}
{\parindent=0pt \v5 Let $K$ be a field, and let $R$ be a finitely
generated free $K$-algebra, or a path algebra defined by a finite
directed graph over $K$. Then it is well-known that $R$ holds an
effective Gr\"obner basis theory ([Mor], [FFG]) which generalizes
successfully the celebrated algorithmic Gr\"obner basis theory for
commutative polynomial algebras [Bu].  Also from the literature
(e.g. [AL], [Ap], [HT], [K-RW], [Li1], [BGV], [Lev]) we know that
many important quotient algebras of $R$, such as exterior algebras,
Clifford algebras, solvable polynomial algebras and some of their
homomorphic images hold an effective Gr\"obner basis theory. In
general for an arbitrary (two-sided) ideal $I$ of $R$, it seems
hardly to know the existence of a Gr\"obner basis theory for the
quotient algebra $A=R/I$. As we learnt from loc. cit., the first
step to have a Gr\"obner basis theory for $A$ is to have a ``nice"
$K$-basis for $A$, for instance, a {\it multiplicative $K$-basis}
$\B$ in the sense that $u,v\in\B$ implies $uv=0$ or $uv\in\B$. In
[Gr2] the following result was obtained for algebras defined by a
2-nomial ideal.{\parindent=0pt\v5

{\bf Theorem} ([Gr2], Theorem 2.3) Suppose that $R$ is a $K$-algebra
with multiplicative $K$-basis ${\cal C}$. Let $I$ be an ideal in $R$
and $\pi$: $R\r R/I$ be the canonical surjection. Let $\pi ({\cal
C})^*=\pi ({\cal C})-\{ 0\}$. Then $\pi ({\cal C})^*$ is a
multiplicative $K$-basis for $R/I$ if and only if $I$ is a 2-nomial
ideal (i.e., $I$ is generated by elements of the form $u-v$ and $w$
where $u,v,w\in {\cal C}$). }\v5

Also as shown in [Gr1, 2], if an algebra $A$ has a multiplicative
$K$-basis $\B$ and if there exists a (two-sided) monomial ordering
$\prec$ on $\B$, then a Gr\"obner basis theory may be developed for
$A$ and such a computational ideal theory can be further applied to
develop a computational module theory (representation theory) of
$A$.  In [Li2], the Gr\"obner basis theory for algebras with a
multiplicative $K$-basis [Gr1, 2] was extended to a Gr\"obner basis
theory for algebras with a {\it skew multiplicative $K$-basis} $\B$
(i.e., $u,v\in\B$ implies $uv=0$ or $uv=\lambda w$ for some nonzero
$\lambda\in K$ and $w\in\B$) and was used to study the structure
theory of  quotient algebras of a $\Gamma$-graded $K$-algebra via
their $\Gamma$-leading and $\B$-leading homogeneous algebras.}\par

Inspired by the work of [Gr1, 2] and [Li2], in this paper we look
for a (left, right, or two-sided) Gr\"obner basis theory for skew
2-nomial algebras and almost skew 2-nomial algebras (see the
definitions in Section 3 and Section 4) at a modest level. More
precisely, in Section 1 we give a quick introduction to the
Gr\"obner basis theory of $K$-algebras with a skew multiplicative
$K$-basis (where $K$ is a field). In Section 2 we consider the class
of $K$-algebras defined by monomials and binomials of certain
special types, and we demonstrate respectively the existence of a
left Gr\"obner basis theory for a subclass of such algebras that may
not necessarily have a right Gr\"obner basis theory, the existence
of a right Gr\"obner basis theory for a subclass of such algebras
that may not necessarily have a left Gr\"obner basis theory, and the
existence of a two-sided Gr\"obner basis theory for a subclass of
such algebras. The first two subclasses of algebras we discussed
include numerous quantum binomial algebras  which provide binomial
solutions to the Yang-Baxter equation ([Laf], [G-I2]), and the third
subclass of algebras generalize the well-known multiparameter
quantized coordinate ring of affine $n$-space $K_Q[z_1,...,z_n]$ (a
precise definition is quoted from [GL] in the beginning of Section
1) which is a typical solvable polynomial algebra in the sense of
[K-RW]. Moreover, the practical examples given in this section also
make the basis for us to discuss the possibility of lifting a (left,
right, or two-sided) Gr\"obner basis theory. Motivated by the
results of Section 2 (mainly Theorem 2.2 -- Corollary 2.5), in
Section 3 we introduce general skew 2-nomial algebras, and, as the
first step of having a (left, right, or two-sided) Gr\"obner basis
theory, the existence of a skew multiplicative $K$-basis for every
skew 2-nomial algebra is established. In Section 4 we introduce the
class of almost skew 2-nomial algebras  so that the working
principle of [Li2], combined with the results of Sections 2 -- 3,
may be applied effectively to the algebras with the associated
graded algebra which is a skew 2-nomial algebra (as shown by Example
(8) of Section 2). The final Section 5 is for summarizing several
open problems related to previous sections.
\par

Throughout this paper, $K$ always denotes a field. By a $K$-algebra
we always mean a  {\it finitely generated} associative $K$-algebra
with identity element $1$, that is, the algebra of the form $\Lambda
=K[a_1,...,a_n]$ with the set of generators $\{ a_1,...,a_n\}$.  For
a given $K$-algebra $\Lambda$, if a ``left ideal" or ``right ideal"
of $\Lambda$ is not specified in the text, ideals are all
``two-sided ideals". Moreover, we write $K^*$ for $K-\{ 0\}$, and
write $\langle S\rangle$ for the two-sided ideal generated by a
subset $S$ in the algebra considered. Also we use $\NZ$ to denote
the set of nonnegative integers.\v5

\section*{1. Grobner Bases with Respect to Skew Multiplicative $K$-basis}
Based on [Gr2] and [Li2], we begin with introducing the Gr\"obner
basis theory of $K$-algebras with a skew multiplicative $K$-basis
(see the definition below).\v5

Let $\Lambda =K[a_1,...,a_n]$ be a $K$-algebra generated by
$\{a_1,...,a_n\}$ and let $\B$ be a $K$-basis of $\Lambda$
consisting of elements of the form
$a_{i_1}^{\alpha_1}a_{i_2}^{\alpha_2}\cdots a_{i_n}^{\alpha_n}$
where $a_{i_j}\in\{ a_1,...,a_n\}$ and $\alpha_j\in\NZ$ (in
principle such a $K$-basis always exists). Denoting elements of $\B$
by using lowercase letters $u$, $v$, $w$, ..., if $\B$ has the
property that
$$u,~v\in\B~\hbox{implies}~\left\{\begin{array}{l} u\cdot v=\lambda w~\hbox{for some}
~\lambda\in K^*,~w\in\B ,\\
\hbox{or}~u\cdot v=0,\end{array}\right.$$ then we call $\B$ a  {\it
skew multiplicative $K$-basis} of $\Lambda$.\v5

Obviously, a skew multiplicative $K$-basis generalizes the notion of
a {\it multiplicative $K$-basis} used in representation theory and
classical Gr\"obner basis theory of associative algebras (i.e., a
$K$-basis $\B$ with the property that $u$, $v\in\B$ implies $uv=0$
or $uv\in\B$). Typical algebras with a skew multiplicative $K$-basis
are ordered semigroup algebras, free algebras, commutative
polynomial algebras, path algebras defined by finite directed
graphs, exterior algebras, and the well-known skew polynomial
algebra $K_Q[z_1,...,z_n]$, or more precisely, the multiparameter
quantized coordinate ring of affine $n$-space  (e.g. see [GL])
defined subject to the relations $z_jz_i=q_{ji}z_iz_j$, $1\le i,j\le
n$,  where $Q =(q_{ij})$ is a multiplicatively antisymmetric
$n\times n$ matrix (i.e., $q_{ii}=1$ and $q_{ij}=q_{ji}^{-1}\in K^*$
for $1\le i,j\le n$).  Section 2 and Section 3 will provide a large
class of algebras with a skew multiplicative $K$-basis.\v5

From now on we let $\Lambda$ be a $K$-algebra with a skew
multiplicative $K$-basis $\B$, as defined above. \v5

Let $\prec$ be a total ordering on $\B$. Adopting the commonly used
terminology in computational algebra, we call an element $u\in\B$ a
{\it monomial}; and for each $f\in\Lambda$, say
$$f=\sum^s_{i=1}\lambda_iu_i,\quad \lambda_i\in K^*,~u_i\in\B,~u_1\prec u_2\prec\cdots\prec u_s ,$$
the {\it leading monomial} of $f$, denoted $\LM (f)$, is defined as
$\LM (f)=u_s$; the {\it leading coefficient} of $f$, denoted $\LC
(f)$, is defined as $\LC (f) =\lambda_s$; and the {\it leading term}
of $f$, denoted $\LT (f)$, is defined as $\LT (f)=\LC (f)\LM
(f)=\lambda_su_s$. If $S\subset\Lambda$, then we write $\LM (S)=\{
\LM (f)~|~f\in S\}$ for the set of leading monomials of $S$. If
furthermore $\prec$ is a well-ordering, i.e., a total ordering
satisfying descending chain condition on $\B$ (in principle, though,
such an ordering exists by the well-known well-ordering theorem from
axiomatic set theory), then  a (left) admissible system for
$\Lambda$ may be introduced as follows. {\parindent=0pt\v5

{\bf 1.1. Definition} (i) An ordering $\prec$ on $\B$ is said to be
a {\it left monomial ordering} if $\prec$ is a well-ordering such
that the following conditions are satisfied:\vskip
6pt\parindent=2.4truecm \re{({\bf LMO1})} If $w,u,v\in\B$ are such
that  $w\prec u$ and $\LM (vw),~\LM (vu)\not\in K^*\cup\{ 0\}$, then
$\LM (vw)\prec \LM (vu)$. \re{({\bf LMO2})} For $w,u\in\B$, if
$u=\LM (vw)$ for some $v\in\B$ with $v\ne 1$ (in the case $1\in\B$),
then $w\prec u$.
\parindent=0pt\par
In this case we call the pair $(\B, \prec )$ a {\it left admissible
system} of $\Lambda$. \par

(ii) An ordering $\prec$ on $\B$ is said to be a {\it two-sided
monomial ordering} if $\prec$ is a well-ordering such that the
following conditions are satisfied:\vskip 6pt\parindent=2.05truecm
\re{({\bf MO1})} If $w,u,v,s\in\B$ are such that $w\prec u$ and $\LM
(vws), ~\LM (vus)\not\in K^*\cup\{ 0\}$, then $\LM (vws)\prec \LM
(vus)$. \re{({\bf MO2})} For $w,u\in\B$, if $u=\LM (vws)$ for some
$v,s\in\B$ with $v\ne 1$ or $s\ne 1$ (in the case $1\in\B$), then
$w\prec u$. \parindent=0pt\par
In this case we call the pair $(\B, \prec )$ a {\it two-sided
admissible system} of $\Lambda$. \v5

{\bf Remark} (i) Note that if $\Lambda$ is the path algebra defined
by a finite directed graph with finitely $n\ge 2$ edges
$e_1,...,e_n$, then the multiplicative  $K$-basis $\B$ of $\Lambda$
does not contain the identity element $1=e_1+\cdots+e_n$. That is
why we specified the case  $1\in\B$.\par

(ii) In the case that $1\in\B$, both (LMO2) and (MO2) turn out that
$1\prec u$ for all $u\in\B$ with $u\ne 1$. This shows that our
definition of a (left, two-sided) monomial ordering is consistent
with that used in the classical commutative and noncommutative
cases. \par

(iii) If $1\in\B$ and $u,v\in \B-\{ 1\}$ satisfy $uv=\lambda\in
K^*$, then $\LM (uv)=1$. Noticing (ii) above, we must ask
$\LM(uv)\not\in K^*\cup\{ 0\}$ in order to make (LMO1) and (MO1)
valid, . }\v5

If the algebra $\Lambda=\oplus_{\gamma\in\Gamma}\Lambda_{\gamma}$ is
a $\Gamma$-graded $K$-algebra, where $\Gamma$ is an ordered
semigroup by a total ordering $<$, such that the skew multiplicative
$K$-basis of $\Lambda$ {\it consists of $\Gamma$-homogeneous
elements}, and if $\prec$ is a well-ordering on $\B$, then we may
define the ordering $\prec_{gr}$ for $u$, $v\in\B$ subject to the
rule:
$$\begin{array}{rcl} u\prec_{gr} v&\Leftrightarrow& d(u)<d(v)\\
&~~&\hbox{or}~d(u)=d(v)~\hbox{and}~u\prec v,\end{array}$$ where
$d(~)$ is referred to as the degree function on homogeneous elements
of $\Lambda$. If $\prec_{gr}$ is a (left or two-sided) monomial
ordering on $\B$ in the sense of Definition 1.1, then we call
$\prec_{gr}$ a (left or two-sided) {\it $\Gamma$-graded monomial
ordering} on $\B$ (or on $\Lambda$). A typical $\NZ$-graded (left or
two-sided) monomial ordering is the $\NZ$-graded (reverse)
lexicographic ordering on a free $K$-algebra $K\langle
X_1,...,X_n\rangle$, a commutative polynomial $K$-algebra
$K[x_1,...,x_n]$, and a path algebra $KQ$ defined by a finite
directed graph $Q$ over $K$, where the $\NZ$-gradation may be any
{\it weight $\NZ$-gradation} obtained by assigning  to each
generator a {\it positive degree} $n_i$, for, in each case
considered the ``monomials" in the standard $K$-basis $\B$ are all
$\NZ$-homogeneous elements.{\parindent=0pt\v5

{\bf Convention} Here let us make the convention once and for all
that a lexicographic ordering always means a ``left lexicographic
ordering" whenever such a ordering is used.}\v5

Since $\B$ is a skew multiplicative $K$-basis of $\Lambda$, the
division of monomials is well-defined, namely, for $u,v\in\B$, we
say that $u$ {\it divides} $v$, denoted $u|v$, provided there are
$w,s\in\B$ and $\lambda\in K^*$ such that $v=\lambda wus$;
similarly, for $u,v\in\B$, we say that {\it $u$ divides $v$ from
left-hand side}, also denoted $u|v$, provided there is some $w\in\B$
and $\lambda\in K^*$ such that $v=\lambda wu$. If furthermore
$\Lambda$ has a (left) admissible system $(\B ,\prec)$, then the
division of monomials may be extended to the division of two
elements in $\Lambda$ as follows. Let $f,g\in\Lambda$. If $\LM
(g)|\LM (f)$, then there are $w,s\in\B$ and $\lambda\in K^*$ such
that $\LM (f)=\lambda w\LM (g)s$. Writing $g=\LC (g)\LM (g)+g'$ with
$\LM (g')\prec \LM (g)$ and $\LC (wgs-wg's)=\mu$, we have
$$f=\frac{\LC (f)}{\mu\LC(g)}wgs+f'~\hbox{satisfying}~\LM (f)=\LM (w\LM (g)s),~\LM (f')\prec\LM (f).$$
If $\LM (g)|\LM (f')$, repeat the same procedure for $f'$. While in
the case that $\LM (g)\not |~\LM (f)$, turn to consider the
divisibility of $\LM (f_1)$ by $\LM (g)$, where $f_1=f-\LT (f)$.
Since $\prec$ is a well-ordering, the division process above stops
after a finite number of steps and eventually we obtain
$$\begin{array}{rcl} f&=&\sum_{i}\lambda_{i}w_{i}gs_{i}+r_f~\hbox{with}~
\lambda_i\in K^*,~w_i,s_i\in\B ,~r_f\in\Lambda ,\\
&{~}&\hbox{satisfying}~\LM (w_igs_i)\preceq\LM (f),~\LM
(r_f)\preceq\LM
(f),\\
&{~}&\hbox{and if}~r_f=\sum_j\mu_jw_j~\hbox{then none of
the}~w_j'\hbox{s can be  divided by}~\LM (g).\end{array}$$\par The
division procedure of $f$ by $g$ from left-hand side can be
performed in a similar way, and the output result is
$$\begin{array}{rcl} f&=&\sum_{i}\lambda_{i}w_{i}g+r_f~\hbox{with}~
\lambda_i\in K^*,~w_i\in\B ,~r_f\in\Lambda ,\\
&{~}&\hbox{satisfying}~\LM (w_ig)\preceq\LM (f),~\LM (r_f)\preceq\LM
(f),\\
&{~}&\hbox{and if}~r_f=\sum_j\mu_jw_j~\hbox{then none of
the}~w_j'\hbox{s can be  divided by}~\LM (g).\end{array}$$\par
Actually the division procedure above gives rise to an effective
division algorithm manipulating the reduction of elements by a
subset of $\Lambda$, and thereby leads to the notion of a (left)
Gr\"obner basis.{\parindent=0pt\v5

{\bf 1.2. Definition}  (i) Let $(\B ,\prec )$ be an admissible
system of the algebra $\Lambda$, and $I$ an ideal of $\Lambda$. A
subset $\G\subset I$ is said to be a {\it Gr\"obner basis} of $I$
(or a {\it two-sided Gr\"obner basis} of $I$ in case the phrase
``two-sided" is need to be emphasized) if for each nonzero $f\in I$,
there is some $g\in\G$ such that $\LM (g)|\LM (f)$, or equivalently,
if each nonzero $f\in I$ has a {\it Gr\"obner presentation} by
elements of $\G$, i.e.,
$$\begin{array}{rcl} f&=&\SUM^{}_{i,j}\lambda_{ij}w_{ij}g_jv_{ij}~\hbox{with}~
\lambda_{ij}\in K^*,~w_{ij},v_{ij}\in\B,~g_j\in \G,~\hbox{satisfying}\\
&{~}&\LM (w_{ij}g_jv_{ij})\preceq\LM (f),~\hbox{and}~\LM
(w_{ij^*}\LM (g_{j^*})v_{ij^*})=\LM (f)~\hbox{for
some}~j^*.\end{array}$$ (ii) Let $(\B ,\prec )$ be a left admissible
system of the algebra $\Lambda$, and $L$ a left ideal of $\Lambda$.
A subset $\G\subset L$ is said to be a {\it left Gr\"obner basis} of
$L$ if for each nonzero $f\in L$, there is some $g\in\G$ such that
$\LM (g)|\LM (f)$ from left-hand side, or equivalently, if each
nonzero $f\in L$ has a {\it left Gr\"obner presentation} by elements
of $\G$, i.e.,
$$\begin{array}{rcl} f&=&\SUM^{}_{i,j}\lambda_{ij}w_{ij}g_j~\hbox{with}~
\lambda_{ij}\in K^*,~w_{ij}\in\B,~g_j\in \G,~\hbox{satisfying}\\
&{~}&\LM (w_{ij}g_j)\preceq\LM (f),~\hbox{and}~\LM (w_{ij^*}\LM
(g_{j^*}))=\LM (f)~\hbox{for some}~j^*.\end{array}$$} \v5

Let $(\B ,\prec )$ be a (left) admissible system of the algebra
$\Lambda$. A subset $\Omega\subset\B$ is said to be {\it reduced} if
$u,v\in\Omega$ and $u\ne v$ implies $u\not |~ v$. If $\G$ is a
(left) Gr\"obner basis of a (left) ideal $J$ in $\Lambda$ such that
$\LM (\G )$ is reduced, then we say that $\G$ is {\it LM}-{\it
reduced}. It is easy to see that a (left) Gr\"obner basis $\G$ of a
(left) ideal $J$ is LM-reduced if and only if $\G$ is a {\it
minimal} (left) Gr\"obner basis in the sense that any proper subset
of $\G$ cannot be a (left) Gr\"obner basis for $J$.\par

Let $(\B ,\prec )$ be a (left) admissible system of the algebra
$\Lambda$.  If $J$ is a nonzero (left) ideal of $\Lambda$, then
$J^*=J-\{ 0\}$ is trivially a (left) Gr\"obner basis of $J$.
Although we do not know if in general there is an effective way to
construct a ``smaller" (left) Gr\"obner basis for $J$, the next
proposition tells us that  each (left) ideal of $\Lambda$ has a {\it
minimal} (left) Gr\"obner basis $\G$.{\parindent=0pt\v5

{\bf 1.3. Proposition} (i) Suppose that the algebra $\Lambda$ has an
admissible system $(\B ,\prec )$. Then each ideal $I$ of $\Lambda$
has a minimal Gr\"obner basis
$$\G =\{ g\in I~|~\hbox{if}~g'\in I~\hbox{and}~g'\ne g,~\hbox{then}~\LM (g')\not |~\LM (g)\} .$$
(ii) Suppose that the algebra $\Lambda$ has a left admissible system
$(\B ,\prec )$. Then each left ideal $L$ of $\Lambda$ has a minimal
left Gr\"obner basis
$$\G =\{ g\in L~|~\hbox{if}~g'\in L~\hbox{and}~g'\ne g,~\hbox{then}~\LM (g')\not |~\LM (g)\} .$$\par\QED

\v5 {\bf 1.4. Definition} Let $\Lambda$ be a $K$-algebra with a skew
multiplicative $K$-basis $\B$. If $\Lambda$ has a left, respectively
a two-sided admissible system $(\B ,\prec )$, then we say that
$\Lambda$ has a {\it left}, respectively a {\it two-sided} {\it
Gr\"obner basis theory}.}\par

If $\Lambda$ has a left,  respectively a two-sided Gr\"obner basis
theory such that every left, respectively every two-sided ideal has
a finite left, respectively a finite two-sided Gr\"obner basis, then
we say that $\Lambda$ has a {\it finite left}, respectively a {\it
finite two-sided} {\it Gr\"obner basis theory}. \v5

Suppose that the algebra $\Lambda$ has a (left, or two-sided)
Gr\"obner basis theory in the sense of Definition 1.4, and let $J$
be a (left, or two-sided) {\it monomial ideal} of $\Lambda$, i.e.,
$J$ is generated by a subset $\Omega\subset\B$. Since $\B$ is a skew
multiplicative $K$-basis, $J$ behaves exactly like a monomial ideal
in a commutative polynomial algebra, that is, $u\in\B\cap J$ if and
only if there exists some $v\in\Omega$ such that $v|u$; and for an
element $f=\sum_i\lambda_iu_i\in\Lambda$ with $\lambda_i\in K^*$,
$u_i\in\B$, $f\in J$ if and only if all $u_i\in J$, the Gr\"obner
basis theory of $\Lambda$ is featured as follows.{\parindent=0pt\v5

{\bf 1.5. Proposition} (i) Suppose that the algebra $\Lambda$ has a
left Gr\"obner basis theory, and let $L$ be a left ideal of
$\Lambda$. A subset $\G\subset L$ is a left Gr\"obner basis of $L$
if and only if ${~}_{_{\Lambda}}\langle\LM (L)\rangle
={~}_{_{\Lambda}}\langle\LM (\G )\rangle$ where both sides of the
equality are left ideals  generated by $\LM (L)$, respectively by
$\LM (\G )$; $L$ has a finite left  Gr\"obner basis basis if and
only if ${~}_{_{\Lambda}}\langle\LM (L)\rangle$ has a finite
monomial generating set. \par

(ii) Suppose that the algebra $\Lambda$ has a two-sided Gr\"obner
basis theory, and let $I$ be an ideal of $\Lambda$. A subset
$\G\subset I$ is a Gr\"obner basis of $I$ if and only if $\langle\LM
(I)\rangle =\langle\LM (\G )\rangle$; $I$ has a finite Gr\"obner
basis if and only if $\langle\LM (I)\rangle$ has a finite monomial
generating set.\par\QED}\v5

Finally, if a $K$-algebra $\Lambda$ has a (left) Gr\"obner basis
theory in the sense of Definition 1.4, then the {\it fundamental
decomposition theorem} holds, that is, if $I$ is a (left) ideal of
$\Lambda$, then the $K$-vector space $\Lambda$ is decomposed into a
direct sum of two subspaces:
$$\Lambda =I\oplus K\hbox{-span}N(I),$$
where $N(I)=\{ u~|~u\in\B ,~\LM (f)\not |~u,~f\in I\}$ is the set of
normal monomials in $\B$ (modulo $I$). In the case that $I$ is
generated by a (left) Gr\"obner basis $\G$, $N(I)$ can be obtained
by means of the (left) division algorithm by $\LM (\G
)$.{\parindent=0pt\v5

{\bf Remark} Let $\Lambda$ be a $K$-algebra with a skew
multiplicative $K$-basis $\B$. To end this section, it is necessary
to note the following:\par

(i) A {\it right Gr\"obner basis theory} for right ideals in
$\Lambda$ can be stated in a similar way.\par

(ii) From Definition 1.1 it is clear that if the $K$-basis $\B$ of
$\Lambda$ contains the identity element 1, then any two-sided
monomial ordering on $\B$ is also a one-sided (i.e. left or right)
monomial ordering, and thereby a two-sided Gr\"obner basis theory
for $\Lambda$ is certainly a one-sided Gr\"obner basis theory for
$\Lambda$. But conversely, in the next section we will see the
examples showing that a one-sided monomial ordering is not
necessarily a two-sided monomial ordering (note that such a proper
example has been missing on page 35 of [Li1]). }\v5

\section*{2. Some Motivating Results}
Let $R=K[a_1,...,a_n]$ be a $K$-algebra generated by $\{
a_1,...,a_n\}$. Assume that $R$ has a two-sided Gr\"obner basis
theory with respect to a two-sided admissible system $(\B ,\prec)$
in the sense of Definition 1.4, where $\B$ is a {\it skew
multiplicative} $K$-{\it basis} of $R$  and $\prec$ is a  {\it
two-sided monomial ordering} on $\B$. For instance, $R$ is a
commutative polynomial algebra, a free algebra, a path algebra, or
the multiparameter quantized coordinate ring of affine $n$-space
$K_Q[z_1,...,z_n]$ with the parameter matrix $Q$ as described in the
beginning of Section 1. In this section, firstly we consider the
quotient algebra $R/I$ of $R$, where $I$ is an ideal generated by
monomials and binomials of certain special types, and we demonstrate
respectively the existence of a left Gr\"obner basis theory for a
subclass of such algebras that may not necessarily have a two-sided
Gr\"obner basis theory, the existence of a right Gr\"obner basis
theory for a subclass of such algebras that may not necessarily have
a two-sided Gr\"obner basis theory, as well as the existence of a
two-sided Gr\"obner basis theory for a subclass of such algebras;
and then, in the case where $R=\oplus_{\gamma\in\Gamma}R_{\gamma}$
is a $\Gamma$-graded algebra by an ordered semigroup $\Gamma$ which
is ordered by a well-ordering, for an arbitrary ideal $I$ of $R$, we
apply the foregoing results to  the $\Gamma$-leading homogeneous
algebra $A^{\Gamma}_{\rm LH}=R/\langle\LH (I)\rangle$ of $A$ (see
the definition given before Proposition 2.6); moreover, we show that
if $A^{\Gamma}_{\rm LH}$ is a domain, then a (left, right, or
two-sided) Gr\"obner basis theory of $A^{\Gamma}_{\rm LH}$ may be
lifted to a (left, right, or two-sided) Gr\"obner basis theory for
$A$. Here and in what follows, with respect to a skew multiplicative
$K$-basis, a one-sided or two-sided Gr\"obner basis theory always
means the one in the sense of Definition 1.4,  otherwise it is the
one indicated before Theorem 2.10.\v5

Let $I$ be an ideal of $R$, and $A=R/I$. Since $R$ has a Gr\"obner
basis theory, under the canonical algebra epimorphism $\pi$: $R\r
A$, the set $N(I)$ of normal monomials in $\B$ (modulo $I$) projects
to a $K$-basis of $A$, that is $\{ \bar{u}=u+I~|~u\in N(I)\}$ (see
Section 1). In what follows we use $\OV{N(I)}$ to denote this basis.
\par Our first result is to deal with monomial ideals, i.e., ideals
generated by elements of $\B$. {\parindent=0pt\v5

{\bf 2.1. Proposition} Let $R$ and $(\B ,\prec )$ be as fixed above.
Consider a subset $\Omega\subset\B$ (where $1\not\in\Omega$ if
$1\in\B$) and the ideal $I=\langle\Omega\rangle$ in $R$. The
following statements hold.\par (i) $\OV{N(I)}$ is a skew
multiplicative $K$-basis for the quotient algebra $A=R/I$.\par (ii)
The two-sided monomial ordering $\prec$ on $\B$ induces a two-sided
monomial ordering on $\OV{N(I)}$, again denoted $\prec$, and hence
the quotient $A=R/I$ has a two-sided  Gr\"obner basis theory.\vskip
6pt {\bf Proof} (i) Since $\B$ is a skew multiplicative $K$-basis of
$R$, it follows from the remark made preceding Proposition 1.6 that
that $u\in\B\cap I$ if and only if there is some $v\in\Omega$ such
that $v|u$, i.e., $u=\lambda wvs$ with $\lambda\in K^*$ and
$w,s\in\B$. So, suppose $\OV{u_1}$, $\OV{u_2}\in\OV{N(I)}$ and
$\OV{u_1}\OV{u_2}\ne 0$, then $u_1u_2=\mu v$ with $\mu\in K^*$ and
$v\in N(I)$. Hence $\OV{u_1}\OV{u_2}=\mu\OV{v}$. This shows that
$\OV{N(I)}$ is a skew multiplicative $K$-basis for $A$.\par

(ii) Suppose $\bar{u},\bar{v}\in\OV{N(I)}$ such that
$\bar{u}\prec\bar{v}$. Then $u\prec v$  since we are using the
induced ordering $\prec$ on $\OV{N(I)}$. If
$\bar{w},\bar{s}\in\OV{N(I)}$ such that $\LM
(\bar{w}\bar{u}\bar{s})\not\in K^*\cup \{ 0\}$ and $\LM
(\bar{w}\bar{v}\bar{s})\not\in K^*\cup\{ 0\}$, then
$\bar{w}\bar{u}\bar{s}\not\in K^*\cup\{ 0\}$,
$\bar{w}\bar{v}\bar{s}\not\in K^*\cup\{ 0\}$, and as argued in the
proof of (i) above, there are $u',v'\in N(I)$ such that $wus=\lambda
u'$ and $wvs=\mu v'$, where $\lambda ,\mu\in K^*$. Hence, it follows
from $\LM (wus)=u'\prec v'=\LM (wvs)$ that $\LM
(\bar{w}\bar{u}\bar{s})=\OV{u'}\prec \OV{v'}=\LM
(\bar{w}\bar{v}\bar{s})$. This shows that (MO1) of Definition 1.1
holds. To prove that (MO2) of Definition 1.1 holds as well, suppose
$\bar{u},\bar{v},\bar{w},\bar{s}\in\OV{N(I)}$ such that $\bar{u}=\LM
(\bar{w}\bar{v}\bar{s})$ (where $\bar{w}\ne 1$ or $\bar{s}\ne 1$ if
$1\in\B$, and thereby $w\ne 1$ or $s\ne 1$ by the choice of
$\Omega$), and $wvs=\lambda v'$ with $\lambda\in K^*$ and $v'\in
N(I)$. Thus, $v'=\LM (wvs)$ implies $v\prec v'$, and consequently,
$\bar{v}\prec \OV{v'}=\LM (\bar{w}\bar{v}\bar{s})=\bar{u}$, as
desired. \QED}\v5

Next we show that there is a left Gr\"obner basis theory for a class
of algebras defined by monomials and binomials of certain special
type. {\parindent=0pt\v5

{\bf 2.2. Theorem} Let $R$ and $(\B ,\prec )$ be as fixed above.
Suppose that the skew multiplicative $K$-basis $\B$ contains every
monomial of the form
$a_{\ell_1}^{\alpha_1}a_{\ell_2}^{\alpha_2}\cdots
a_{\ell_n}^{\alpha_n}$ with respect to some permutation $a_{\ell_1},
a_{\ell_2}, ..., a_{\ell_n}$ of $a_1, a_2, ..., a_n$, where
$\alpha_1,...,\alpha_n\in\NZ$. If $\G =\Omega\cup G$ is a minimal
Gr\"obner basis of the ideal $I=\langle\G \rangle$, where
$\Omega\subset\B$ and $G$ consists of $\frac{n(n-1)}{2}$ elements
$$g_{ji}=a_{\ell_j}a_{\ell_i}-\lambda_{ji}a_{\ell_i}a_{\ell_p}, ~
1\le i<j\le n,~i<p\le n,~ \lambda_{ji}\in K^*\cup \{ 0\},$$ such
that $\LM (g_{ji})=a_{\ell_j}a_{\ell_i}$, $1\le i<j\le n$, then  the
following statements hold for the algebra $A=R/I$.\par

(i) $N(I)\subseteq \{
a_{\ell_1}^{\alpha_1}a_{\ell_2}^{\alpha_2}\cdots
a_{\ell_n}^{\alpha_n}~|~\alpha_1,...,\alpha_n\in\NZ\}$, and
$\OV{N(I)}$ is a skew multiplicative $K$-basis for $A$.\par

(ii) Let $\prec_{_I}$ denote the ordering on $\OV{N(I)}$, which is
induced by the lexicographic ordering or the $\NZ$-graded
lexicographic ordering on $\NZ^n$, such that
$$\OV{a_{\ell_n}}\prec_{_I}\OV{a_{\ell_{n-1}}}\prec_{_I}\cdots
\prec_{_I}\OV{a_{\ell_2}}\prec_{_I}\OV{a_{\ell_1}}.$$
If $\OV{a_{\ell_1}}^{\alpha_1}\OV{a_{\ell_2}}^{\alpha_2}\cdots
\OV{a_{\ell_n}}^{\alpha_n}\ne 0$ implies
$a_{\ell_1}^{\alpha_1}a_{\ell_2}^{\alpha_2}\cdots
a_{\ell_n}^{\alpha_n}\in N(I)$, then $\prec_{_I}$ is a left monomial
ordering on $\OV{N(I)}$, and hence $A$ has a left Gr\"obner basis
theory with respect to the left admissible system
$(\OV{N(I)},\prec_{_I} )$. \vskip 6pt

{\bf Proof} Note that for the convenience of statement, we may
assume, without loss of generality, that $a_{\ell_1}=a_1$,
$a_{\ell_2}=a_2$, ..., $a_{\ell_n}=a_n$.\par

(i) By the assumptions, $\B$ is a skew multiplicative $K$-basis of
$R$, $a_1,...,a_n\in\B$, and $\LM (g_{ji})=a_ja_i$, $1\le i<j\le n$.
It follows from the division algorithm by $\G$ that $N(I)\subseteq
\{ a_{1}^{\alpha_1}a_{2}^{\alpha_2}\cdots
a_{n}^{\alpha_n}~|~\alpha_j\in\NZ\}$. For $\bar{u}$,
$\bar{v}\in\OV{N(I)}$, if $\bar{u}\bar{v}\ne 0$, then $uv\ne 0$ and
$uv=\lambda w$ for some $\lambda\in K^*$ and $w\in\B$. If $w\in
N(I)$, then $\bar{u}\bar{v}=\lambda\bar{w}$; otherwise, noticing
that $\G$ consists of monomials and binomials, the division of the
monomial $w$ by $\G$ yields a unique Gr\"obner presentation (see
Definition 1.2):
$$w=\sum_{i,j}\mu_{ij}s_{ij}g_{ji}t_{ij}+\eta u',~\hbox{where}~
\mu_{ij},\eta\in K^*,~s_{ij},t_{ij}\in\B ,~u'\in N(I),$$ and
consequently $\bar{u}\bar{v}=\lambda\bar{w}=(\lambda\eta)\OV{u'}$.
Therefore, $\OV{N(I)}$ is a skew multiplicative $K$-basis for the
quotient algebra $A=R/I$.\par

(ii) Let $\prec_{_I}$ be one of the orderings on $\OV{N(I)}$ as
mentioned in the theorem, that is, $\prec_{_I}$ is defined subject
to the lexicographic ordering
$$\OV{a_n}\prec_{_I}\OV{a_{n-1}}\prec_{_I}\cdots\prec_{_I}\OV{a_2}\prec_{_I}\OV{a_1}$$
If $\alpha =(\alpha_1,...,\alpha_n)\in\NZ^n$, then we write $|\alpha
|=\alpha_1+\cdots +\alpha_n$.}\par

First note that if $\OV{a_k}\in\{ \OV{a_1}, \OV{a_2},...,\OV{a_n}\}$
and if $\OV{a_k}\cdot\OV{a_i}^{\alpha_i}\not\in K^*\cup\{0\}$,
 then
$$\LM (\OV{a_k}\cdot\OV{a_i}^{\alpha_i})=\left\{\begin{array}{ll}
\OV{a_i}^{\alpha_i+1},&k=i\\ \OV{a_k}\cdot\OV{a_i}^{\alpha_i},&k<i\\
\OV{a_i}^{\alpha_i}\cdot\OV{a_{p_{_k}}}~\hbox{with}~i<p_{_k},&k>i\end{array}\right.
\leqno{(1)}$$
Let $\bar{u}=\OV{a_1}^{\alpha_1}\cdots \OV{a_n}^{\alpha_n}$,
$\bar{v}=\OV{a_1}^{\beta_1}\cdots \OV{a_n}^{\beta_n}\in\OV{N(I)}$ be
such that $\bar{u}\prec_{_I}\bar{v}$. Then
$$\begin{array}{l} \hbox{either}~|\alpha |<|\beta |\\
\hbox{or}~ |\alpha |=|\beta
|~\hbox{and}~\alpha_1=\beta_1,~...,~\alpha_{i-1}=
\beta_{i-1}~\hbox{but}~\alpha_i<\beta_i~\hbox{for
some}~i.\end{array}\leqno{(2)}$$ Now, if
$\bar{s}=\OV{a_1}^{\gamma_1}\cdots\OV{a_n}^{\gamma_n}\in\OV{N(I)}$,
and if $\bar{s}\cdot\bar{u}$, $\bar{s}\cdot\bar{v}\not\in K^*\cup\{
0\}$, then, by (i) we may write
$\LM(\bar{s}\cdot\bar{u})=\OV{a_1}^{\eta_1}\cdots\OV{a_n}^{\eta_n}$,
$\LM
(\bar{s}\cdot\bar{v})=\OV{a_1}^{\rho_1}\cdots\OV{a_n}^{\rho_n}$. It
follows from the foregoing formula (1) that $|\eta |=|\alpha
|+|\gamma |$, $|\rho |=|\beta |+|\gamma |$. Thus, $|\alpha
|<|\beta|$ implies $|\eta |<|\rho |$, and hence
$$\LM(\bar{s}\cdot\bar{u})\prec_{_I}\LM (\bar{s}\cdot\bar{v})\leqno{(3)}$$
provided we are using the graded lexicographic ordering. In the case
that $|\alpha |=|\beta |$, we have $\alpha_1=\beta_1$, ...,
$\alpha_{i-1}=\beta_{i-1}$ and $\alpha_i<\beta_i$ for some $i$ by
the indication (2) given above. Since $\bar{s}\cdot\bar{u}\ne 0$ and
$\bar{s}\cdot\bar{v}\ne 0$,
$\OV{a_1}^{\alpha_1}\cdots\OV{a_{i-1}}^{\alpha_{i-1}}\in\OV{N(I)}$
by the assumption of (ii). By (i), if we put $\LM
(\bar{s}\cdot\OV{a_1}^{\alpha_1}\cdots\OV{a_{i-1}}^{\alpha_{i-1}})=\OV{a_1}^{\tau_1}\cdots
\OV{a_n}^{\tau_n}$, then
$$\begin{array}{l} \LM (\bar{s}\cdot\bar{u})=\LM (\OV{a_1}^{\tau_1}\cdots
\OV{a_n}^{\tau_n}\cdot\OV{a_i}^{\alpha_i}\cdots\OV{a_n}^{\alpha_n}),\\
\LM (\bar{s}\cdot\bar{v})=\LM (\OV{a_1}^{\tau_1}\cdots
\OV{a_n}^{\tau_n}\cdot\OV{a_i}^{\beta_i}\cdots\OV{a_n}^{\beta_n}).\end{array}$$
Applying the previously derived formula (1) to $i+1\le k\le n$, we
have
$\OV{a_k}^{\tau_k}\cdot\OV{a_i}^{\alpha_i}=\OV{a_i}^{\alpha_i}\cdot\OV{a_{p_{_k}}}^{\tau_k}$,
 and
$\OV{a_k}^{\tau_k}\cdot\OV{a_i}^{\beta_i}=\OV{a_i}^{\beta_i}\cdot\OV{a_{p_{_k}}}^{\tau_k}$
with $i<p_k$. This turns out that
$$\begin{array}{l} \LM (\bar{s}\cdot\bar{u})=\LM (\OV{a_1}^{\tau_1}\cdots
\OV{a_i}^{\tau_i}\cdot\OV{a_i}^{\alpha_i}\cdot\OV{a_{p_{i+1}}}^{\tau_{i+1}}
\cdots\OV{a_{p_{_n}}}^{\tau_n}\cdot\OV{a_{i+1}}^{\alpha_{i+1}}\cdots\OV{a_n}^{\alpha_n}),\\
\LM (\bar{s}\cdot\bar{v})=\LM (\OV{a_1}^{\tau_1}\cdots
\OV{a_i}^{\tau_i}\cdot\OV{a_i}^{\beta_i}\cdot\OV{a_{p_{i+1}}}^{\tau_{i+1}}
\cdots\OV{a_{p_{_n}}}^{\tau_n}\cdot\OV{a_{i+1}}^{\beta_{i+1}}\cdots\OV{a_n}^{\beta_n}).\end{array}$$
Hence we obtain
$$\LM(\bar{s}\cdot\bar{u})\prec_{_I}\LM (\bar{s}\cdot\bar{v})\leqno{(4)}$$
This shows that (LMO1) of Definition 1.1(ii) holds. \par

Next, let $\bar{u}=\OV{a_1}^{\alpha_1}\cdots\OV{a_n}^{\alpha_n}$,
$\bar{v}=\OV{a_1}^{\beta_1}\cdots\OV{a_n}^{\beta_n}$,
$\bar{w}=\OV{a_1}^{\gamma_1}\cdots\OV{a_n}^{\gamma_n}\in\OV{N(I)}$
be such that $\bar{u}=\LM (\bar{v}\cdot\bar{w})$ (where $\bar{v}\ne
1$ if $1\in\B$). If $\bar{u}\preceq_{_I}\bar{w}$, then since (LMO1)
holds, we would have $\LM (\bar{v}\cdot\bar{u})\preceq_{_I}\LM
(\bar{v}\cdot\bar{w})=\bar{u}$. But then, from the argument given
before the formula (3) above we would have $|\beta |+|\alpha |\le
|\beta |+|\gamma |=|\alpha |$, that is clearly impossible, for,
$\bar{v}\ne 0$ (also $\bar{v}\ne 1$ if $1\in\B$) implies $|\beta
|\ne 0$. So we must have $\bar{w}\prec_{_I}\bar{u}$. It follows that
(LMO2) of Definition 1.1(ii) holds as well. Summing up,  we conclude
that $\prec_{_I}$ is a left monomial ordering on $\OV{N(I)}$, and
thereby $A$ has a left Gr\"obner basis theory with respect to the
left admissible system $(\OV{N(I)},\prec_{_I})$. \QED\v5

In a similar way, we may obtain a right Gr\"obner basis theory for a
class of algebras defined by monomials and binomials of certain
special type.{\parindent=0pt\v5

{\bf 2.3. Theorem} Under the assumptions on $R$ and $\B$ as in
Theorem 2.2, if $\G=\Omega\cup G$ is a minimal Gr\"obner basis of
the ideal $I=\langle\G \rangle$, where $\Omega\subset\B$ and $G$
consists of $\frac{n(n-1)}{2}$ elements
$$g_{ji}=a_{\ell_j}a_{\ell_i}-\lambda_{ji}a_{\ell_p}a_{\ell_j},
~1\le i<j\le n,~p<j,~ \lambda_{ji}\in K^*\cup \{ 0\},$$ such that
$\LM (g_{ji})=a_{\ell_j}a_{\ell_i}$, $1\le i<j\le n$, then the
following statements hold for the algebra $A=R/I$.\par

(i) $N(I)\subseteq \{
a_{\ell_1}^{\alpha_1}a_{\ell_2}^{\alpha_2}\cdots
a_{\ell_n}^{\alpha_n}~|~\alpha_1,...,\alpha_n\in\NZ\}$, and
$\OV{N(I)}$ is a skew multiplicative $K$-basis for $A$.\par

(ii) Let $\prec_{_I}$ denote the ordering on $\OV{N(I)}$, which is
induced by the $\NZ$-graded reverse lexicographic ordering on
$\NZ^n$, such that
$$\OV{a_{\ell_n}}\prec_{_I}\OV{a_{\ell_{n-1}}}\prec_{_I}\cdots
\prec_{_I}\OV{a_{\ell_2}}\prec_{_I}\OV{a_{\ell_1}}.$$
If $\OV{a_{\ell_1}}^{\alpha_1}\OV{a_{\ell_2}}^{\alpha_2}\cdots
\OV{a_{\ell_n}}^{\alpha_n}\ne 0$ implies
$a_{\ell_1}^{\alpha_1}a_{\ell_2}^{\alpha_2}\cdots
a_{\ell_n}^{\alpha_n}\in N(I)$, then $\prec_{_I}$ is a right
monomial ordering on $\OV{N(I)}$, and hence $A$ has a right
Gr\"obner basis theory with respect to the right admissible system
$(\OV{N(I)},\prec_{_I} )$.\par\QED}\v5

Below we give examples to show that in general algebras of the type
as described in Theorem 2.2 may not necessarily have a two-sided
monomial ordering on $\OV{N(I)}$, and thereby each left monomial
ordering $\prec_{_I}$ used in Theorem 2.2 is in general not a right
monomial ordering on $\OV{N(I)}$;  and that in general algebras of
the type as described in Theorem 2.3 may not necessarily have a
two-sided monomial ordering on $\OV{N(I)}$, and thereby the right
monomial ordering $\prec_{_I}$ used in Theorem 2.3 is in general not
a left monomial ordering on $\OV{N(I)}$. {\parindent=0pt\v5

{\bf Example} (1) In the free $K$-algebra $\KS=K\langle
X_1,X_2,X_3\rangle$ consider the subset $\G$ consisting of
$$\begin{array}{l} g_{21}=X_2X_1-\lambda X_1X_3\\
g_{31}=X_3X_1-\mu X_1X_2\\
g_{32}=X_3X_2-\gamma X_2X_3\end{array}\quad\hbox{where}~\lambda ,\mu
,\gamma\in K^*\cup\{ 0\}.$$ Then, with respect to the natural
$\NZ$-graded lexicographic ordering subject to
$X_1\prec_{grlex}X_2\prec_{grlex}X_3$, where each $X_i$ is assigned
to the degree 1, it is straightforward to verify that $\LM
(g_{ji})=X_jX_i$, $1\le i<j\le 3$, $\G$ forms a Gr\"obner basis of
the ideal $I=\langle\G\rangle$ in each of the following cases:
 $$\begin{array}{l}
 (\lambda ,\mu ,\gamma )=(0 ,\mu ,\gamma ) ~\hbox{with arbitrarily chosen}~
 \mu ,\gamma;\\
 (\lambda ,\mu ,\gamma )=(\lambda ,0 ,\gamma ) ~
 \hbox{with arbitrarily chosen}~\lambda ,\gamma;\\
 (\lambda ,\mu ,\gamma )=(\lambda ,0,0) ~\hbox{with arbitrarily chosen}~\lambda;\\
 (\lambda ,\mu ,\gamma )=(0 ,\mu ,0) ~\hbox{with arbitrarily chosen}~\mu;\\
 (\lambda ,\mu ,\gamma )=(\lambda ,\mu ,\gamma )~\hbox{with}~\lambda ,\mu ,\gamma\in K^*,~\gamma^2=1,\end{array}$$
and $N(I)=\{
X_1^{\alpha_1}X_2^{\alpha_2}X_3^{\alpha_3}~|~\alpha_j\in\NZ\}$.
Since $\G$ is of the type required by Theorem 2.2, it follows that
the algebra $A=\KS /I$ has a left Gr\"obner basis theory with
respect to the left admissible system $(\OV{N(I)},\prec_{_I})$,
where $\prec_{_I}$ is induced by the  lexicographic ordering or the
$\NZ$-graded  lexicographic ordering on $\NZ^n$, such that
$\OV{X_3}\prec_{_I}\OV{X_2}\prec_{_I}\OV{X_1}$. But since the
algebra $A$ has $1,\OV{X_1}$, $\OV{X_2}$, $\OV{X_3}\in\OV{N(I)}$
such that
$$\begin{array}{l} \OV{X_2}\cdot\OV{X_1}=\lambda\OV{X_1}\cdot\OV{X_3},\\
\OV{X_3}\cdot\OV{X_1}=\mu\OV{X_1}\cdot\OV{X_2},\end{array}$$ we see
that if $\lambda\ne 0$, $\mu\ne 0$, then any total ordering $<$ such
that $\OV{X_3}<\OV{X_2}$ or $\OV {X_2}<\OV{X_3}$ cannot be a
two-sided monomial ordering on $\OV{N(I)}$. Therefore, $\prec_{_I}$
cannot be a right monomial ordering on $\OV{N(I)}$.\v5

{\bf Example} (2) In the free $K$-algebra $\KS=K\langle
X_1,X_2,X_3\rangle$ consider the subset $\G$ consisting of
$$\begin{array}{l} g_{21}=X_2X_1-\lambda X_1X_2\\
g_{31}=X_3X_1-\mu X_2X_3\\
g_{32}=X_3X_2-\gamma X_1X_3\end{array}\quad\hbox{where}~\lambda ,\mu
,\gamma\in K^*\cup\{ 0\}.$$ Then, with respect to the natural
$\NZ$-graded lexicographic ordering subject to
$X_1\prec_{grlex}X_2\prec_{grlex}X_3$, where each  $X_i$ is assigned
to the degree 1, it is straightforward to verify that $\LM
(g_{ji})=X_jX_i$, $1\le i<j\le 3$, $\G$ forms a Gr\"obner basis of
the ideal $I=\langle\G\rangle$  in each of the following cases:
 $$\begin{array}{l}
 (\lambda ,\mu ,\gamma )=(0 ,0,\gamma )~\hbox{with arbitrarily chosen}~\gamma;\\
 (\lambda ,\mu ,\gamma )=(0 ,\mu ,0)~\hbox{with arbitrarily chosen}~\mu;\\
 (\lambda ,\mu ,\gamma )=(\lambda ,0 ,\gamma )~\hbox{with arbitrarily chosen}~\lambda ,\gamma;\\
 (\lambda ,\mu ,\gamma )=(\lambda ,\mu ,0 )~\hbox{with arbitrarily chosen}~\lambda ,\mu;\\
 (\lambda ,\mu ,\gamma )=(\lambda ,\mu ,\gamma )~\hbox{with}~\lambda ,\mu ,\gamma\in K^*,~\lambda^2=1,\end{array}$$
and $N(I)=\{
X_1^{\alpha_1}X_2^{\alpha_2}X_3^{\alpha_3}~|~\alpha_j\in\NZ\}$.
Since $\G$ is of the type required by Theorem 2.3, it follows that
the algebra $A=\KS /I$ has a right Gr\"obner basis theory with
respect to the right admissible system $(\OV{N(I)},\prec_{_I})$,
where $\prec_{_I}$ is induced by the $\NZ$-graded reverse
lexicographic ordering on $\NZ^n$, such that
$\OV{X_3}\prec_{_I}\OV{X_2}\prec_{_I}\OV{X_1}$. But since the
algebra $A$ has  $1,\OV{X_1}$, $\OV{X_2}$, $\OV{X_3}\in\OV{N(I)}$
such that
$$\begin{array}{l} \OV{X_3}\cdot\OV{X_1}=\mu\OV{X_2}\cdot\OV{X_3}, \\
\OV{X_3}\cdot\OV{X_2}=\gamma \OV{X_1}\cdot\OV{X_3},\end{array}$$ we
see that if $\mu\ne 0$, $\gamma\ne 0$, then any total ordering $<$
such that $\OV{X_1}<\OV{X_2}$ or $\OV {X_2}<\OV{X_1}$ cannot be a
two-sided monomial ordering on $\OV{N(I)}$. Therefore, $\prec_{_I}$
cannot be a left monomial ordering on $\OV{N(I)}$.}\v5

Also we give an example to show that if the ideal
$I=\langle\G\rangle$ has the Gr\"obner basis $\G$ consisting of
binomials but is neither the type as described in Theorem 2.2 nor
the type as described in Theorem 2.3, then the algebra defined by
$\G$ may not have a two-sided  monomial ordering, and each of the
(left or right) monomial orderings used in both theorems cannot be a
left or right monomial ordering for such an algebra.
{\parindent=0pt\v5

{\bf Example} (3) In the free $K$-algebra $\KS=K\langle X_1,X_2,X_3,
X_4\rangle$ consider the subset $\G$ consisting of
$$\begin{array}{ll} g_{21}=X_2X_1- X_1X_2,&g_{41}=X_4X_1- X_2X_3,\\
g_{31}=X_3X_1- X_2X_4,&g_{42}=X_4X_2- X_1X_3,\\
g_{32}=X_3X_2- X_1X_4,&g_{43}=X_4X_3-X_3X_4.\end{array}$$ Then, with
respect to the natural $\NZ$-graded lexicographic ordering subject
to $X_1\prec_{grlex}X_2\prec_{grlex}X_3\prec_{grlex}X_4$, where each
$X_i$ is assigned to the degree 1, it is straightforward to verify
that $\LM (g_{ji})=X_jX_i$, $1\le i<j\le 4$, $\G$ forms a Gr\"obner
basis of the ideal $I=\langle\G\rangle$, and $N(I)=\{
X_1^{\alpha_1}X_2^{\alpha_2}X_3^{\alpha_3}X_4^{\alpha_4}~|~\alpha_j\in\NZ\}$.
Since the algebra $A=\KS /I$ has $1, \OV{X_1}$, $\OV{X_2}$
$\OV{X_3}$, $\OV{X_4}\in\OV{N(I)}$ such that
$$\begin{array}{ll}
\left\{\begin{array}{l}\OV{X_3}\cdot\OV{X_1}=\OV{X_2}\cdot\OV{X_4}\\
 \OV{X_3}\cdot
\OV{X_2}=\OV{X_1}\cdot\OV{X_4}\end{array}\right.
&\left\{\begin{array}{l}
\OV{X_4}\cdot\OV{X_1}=\OV{X_2}\cdot\OV{X_3}\\
\OV{X_4}\cdot\OV{X_2}=\OV{X_1}\cdot\OV{X_3}\end{array}\right.\\
\\
\left\{\begin{array}{l} \OV{X_4}\cdot
\OV{X_1}=\OV{X_2}\cdot\OV{X_3}\\
\OV{X_3}\cdot\OV{X_1}=\OV{X_2}\cdot\OV{X_4}\end{array}\right.
&\left\{\begin{array}{l} \OV{X_3}\cdot
\OV{X_2}=\OV{X_1}\cdot\OV{X_4}\\
\OV{X_4}\cdot\OV{X_2}=\OV{X_1}\cdot\OV{X_3}\end{array}\right.\end{array}$$
it is clear that there is no two-sided monomial ordering on
$\OV{N(I)}$, and that each of the (left or right) monomial orderings
used in Theorem 2.2 and Theorem 2.3 cannot be a (left, right, or
two-sided) monomial ordering on $\OV{N(I)}$.}\v5

In consideration of the examples given above, algebras defined by
monomials and  binomials of the types as described in Theorem 2.2
and Theorem 2.3, that may have a two-sided Gr\"obner basis theory
are necessarily the type similar to the multiparameter quantized
coordinate ring of affine $n$-space  $K_Q[z_1,...,z_n]$ with the
parameter matrix $Q$ as described in the beginning of Section 1.
This leads to the next result. {\parindent=0pt\v5

{\bf 2.4. Theorem} Under the assumptions on $R$ and $\B$ as in
Theorem 2.2, if $\G =\Omega\cup G$ is a minimal Gr\"obner basis of
the ideal $I=\langle\G \rangle$, where $\Omega\subset\B$ and $G$
consists of $\frac{n(n-1)}{2}$ elements
$$g_{ji}=a_{\ell_j}a_{\ell_i}-\lambda_{ji}a_{\ell_i}a_{\ell_j},
~1\le i<j\le n,~ \lambda_{ji}\in K^*\cup \{ 0\},$$ such that $\LM
(g_{ji})=a_{\ell_j}a_{\ell_i}$, $1\le i<j\le n$, then the following
statements hold for the algebra $A=R/I$.\par

(i) $N(I)\subseteq \{
a_{\ell_1}^{\alpha_1}a_{\ell_2}^{\alpha_2}\cdots
a_{\ell_n}^{\alpha_n}~|~\alpha_1,...,\alpha_n\in\NZ\}$, and
$\OV{N(I)}$ is a skew multiplicative $K$-basis for $A$.\par

(ii) If $\OV{a_{\ell_1}}^{\alpha_1}\OV{a_{\ell_2}}^{\alpha_2}\cdots
\OV{a_{\ell_n}}^{\alpha_n}\ne 0$ implies
$a_{\ell_1}^{\alpha_1}a_{\ell_2}^{\alpha_2}\cdots
a_{\ell_n}^{\alpha_n}\in N(I)$, then $A$ has a two-sided Gr\"obner
basis theory with respect to the two-sided admissible system
$(\OV{N(I)},\prec_{_I} )$, where $\prec_{_I}$ is the two-sided
monomial ordering on $\OV{N(I)}$ induced by  one of the following
orderings on $\NZ^n$: the lexicographic ordering, the $\NZ$-graded
lexicographic ordering, and the $\NZ$-graded reverse lexicographic
ordering,  such that
$$\OV{a_{\ell_n}}\prec_{_I}\OV{a_{\ell_{n-1}}}\prec_{_I}\cdots
\prec_{_I}\OV{a_{\ell_2}}\prec_{_I}\OV{a_{\ell_1}}.$$   \vskip 6pt

{\bf Proof} Note that in this case we have
$g_{ji}=a_ja_i-\lambda_{ji}a_ia_j$, $1\le i<j\le n$,
$\lambda_{ji}\in K^*\cup\{ 0\}$. Let $\prec_{_I}$ be one of the
orderings on $\OV{N(I)}$ as mentioned.  If
$\bar{u}=\OV{a_1}^{\alpha_1}\cdots \OV{a_n}^{\alpha_n}$,
$\bar{v}=\OV{a_1}^{\beta_1}\cdots \OV{a_n}^{\beta_n}$,
$\bar{w}=\OV{a_1}^{\gamma_1}\cdots \OV{a_n}^{\gamma_n}$,
$\bar{s}=\OV{a_1}^{\eta_1}\cdots \OV{a_n}^{\eta_n}\in\OV{N(I)}$
satisfy $\bar{u}\prec_{_I}\bar{v}$, $\bar{w}\bar{u}\bar{s}\not\in
K^*\cup \{ 0\}$, and $\bar{w}\bar{v}\bar{s}\not\in K^*\cup\{ 0\}$,
then,  the assumption that $\OV{a_1}^{\alpha_1}\cdots
\OV{a_n}^{\alpha_n}\ne 0$ implies $a_1^{\alpha_1}\cdots
a_n^{\alpha_n}\in N(I)$ and the feature of $g_{ji}$ entail that
$$\LM (\bar{w}\bar{u}\bar{s})=\OV{a_1}^{\gamma_1+\alpha_1+\eta_1}\cdots
\OV{a_n}^{\gamma_n+\alpha_n+\eta_n}\prec_{_I}
\OV{a_1}^{\gamma_1+\beta_1+\eta_1}\cdots
\OV{a_n}^{\gamma_n+\beta_n+\eta_n}=\LM (\bar{w}\bar{v}\bar{s}).$$
This shows that (MO1) of Definition 1.1(i) holds. Similarly,  (MO2)
of Definition 1.1(i) holds as well. Hence $\prec_{_I}$ is a
two-sided monomial ordering on $\OV{N(I)}$, and thereby $A$ has a
two-sided Gr\"obner basis theory with respect to the admissible
system $(\OV{N(I)},\prec_{_I})$. \QED}\v5

Focusing on quotient algebras of free $K$-algebras we have the
following immediate corollary.{\parindent=0pt\v5

{\bf 2.5. Corollary} Let $R=K\langle X_1,...,X_n\rangle$ be the free
$K$-algebra of $n$ generators, $\B$ the standard $K$-basis of $R$,
and let $X_{\ell_1},X_{\ell_2}, ..., X_{\ell_n}$ be a permutation of
$X_1, X_2, ..., X_n$.\par (i) If, with respect to some two-sided
monomial ordering $\prec$ on $\B$, $\G=\Omega\cup G$ is a minimal
Gr\"obner basis of the ideal $I=\langle\G \rangle$, where
$\Omega\subset\B$ and $G$ consists of $\frac{n(n-1)}{2}$ elements
$$g_{ji}=X_{\ell_j}X_{\ell_i}-\lambda_{ji}X_{\ell_i}X_{\ell_p},~1\le i<j\le n,~
i<p\le n,~\lambda_{ji}\in K^*\cup \{ 0\},$$ such that $\LM
(g_{ji})=X_{\ell_j}X_{\ell_i}$, $1\le i<j\le n$, then the conditions
of Theorem 2.2 are satisfied, and consequently the algebra $A=R/I$
 has a left  Gr\"obner basis theory with
respect to the left admissible system $(\OV{N(I)},\prec_{_I})$,
where $\prec_{_I}$ is as described in Theorem 2.2.\par

(ii) If, with respect to some two-sided monomial ordering $\prec$ on
$\B$, $\G =\Omega\cup G$ is a minimal Gr\"obner basis of the ideal
$I=\langle\G\rangle$, where $\Omega\subset\B$ and $G$ consists of
$\frac{n(n-1)}{2}$ elements
$$g_{ji}=X_{\ell_j}X_{\ell_i}-\lambda_{ji}X_{\ell_p}X_{\ell_j},~1\le i<j\le n,~p<j,~
\lambda_{ji}\in K^*\cup \{ 0\},$$ such that $\LM
(g_{ji})=X_{\ell_j}X_{\ell_i}$, $1\le i<j\le n$, then the conditions
of Theorem 2.3 are satisfied, and consequently the algebra $A=R/I$
has a right Gr\"obner basis theory with respect to the right
admissible system $(\OV{N(I)},\prec_{_I})$, where $\prec_{_I}$ is as
described in Theorem 2.3.\par

(iii) If, with respect to some two-sided monomial ordering $\prec$
on $\B$, $\G =\Omega\cup G$ is a minimal Gr\"obner basis of the
ideal $I=\langle\G\rangle$, where $\Omega\subset\B$ and $G$ consists
of $\frac{n(n-1)}{2}$ elements
$$g_{ji}=X_{\ell_j}X_{\ell_i}-\lambda_{ji}X_{\ell_i}X_{\ell_j},~1\le i<j\le n,~
\lambda_{ji}\in K^*\cup \{ 0\},$$ such that $\LM
(g_{ji})=X_{\ell_j}X_{\ell_i}$, $1\le i<j\le n$, then the conditions
of Theorem 2.4 are satisfied, and consequently the  algebra $A=R/I$
has a two-sided Gr\"obner basis theory with respect to the two-sided
admissible system $(\OV{N(I)},\prec_{_I})$, where
 $\prec_{_I}$ is as described in Theorem 2.4.
\par\QED} \v5

The next example illustrates Theorem 2.4, or more precisely,
Corollary 2.5(iii).{\parindent=0pt\v5

{\bf Example} (4)  Consider in the free $K$-algebra $\KS =K\langle
X_1,...,X_n\rangle$ the subset $\G =\Omega\cup G$ consisting of
$$\begin{array}{l} \Omega\subseteq\left\{\left. g_i=X_i^p~ \right |~1\le i\le n\right\} ,
~\hbox{where}~p\ge 2~\hbox{is a fixed integer},\\
G=\left\{\left. g_{ji}=X_jX_i-\lambda_{ji}X_iX_j~\right |~1\le
i<j\le n,~\lambda_{ji}\in K^*\cup\{ 0\}\right\} ,\end{array}$$ and
let $\prec_{grlex}$ be the natural $\NZ$-graded lexicographic
ordering subject to
$X_1\prec_{grlex}X_2\prec_{grlex}\cdots\prec_{grlex}X_n$, where each
$X_i$ is assigned to the degree 1.  Then
$$\LM (\G )=\left\{\left. X_s^p=g_s~\right |~g_s\in\Omega\right\}\cup\left\{\left.
X_jX_i~\right |~1\le i<j\le n\right\} ,$$ and all possible
nontrivial overlap elements of $\G$ are
$$\begin{array}{ll}
S_{i,i\ell}=g_iX_{\ell}-X_i^{p-1}g_{i\ell}=\lambda_{i\ell}X_i^{p-1}X_{\ell}X_i,&i>\ell ,\\
S_{ji,i}=g_{ji}X_i^{p-1}-X_jg_i=-\lambda_{ji}X_iX_jX_i^{p-1},&j>i,\\
S_{j\ell ,\ell t}=g_{j\ell}X_t-X_jg_{\ell t}=\lambda_{\ell
t}X_jX_tX_{\ell}-\lambda_{j\ell}X_{\ell}X_jX_t,&j>\ell
>t.\end{array}$$
For arbitrarily chosen $\lambda_{ji}\in K^*\cup\{ 0\}$, it is
straightforward to verify that the division of $S_{i,i\ell}$,
$S_{ji,i}$, and $S_{j\ell ,\ell t}$ by $\G$ yields respectively the
Gr\"obner presentations (Definition 1.2):
$$\begin{array}{l} S_{i,i\ell}=\SUM^p_{k=2}\lambda_{i\ell}^{k-1}X_i^{p-k}g_{i\ell}X_i^{k-1}+\lambda_{i\ell}^pX_{\ell}g_i,\\
S_{ji,i}=-\SUM^p_{k=2}\lambda_{ji}^{k-1}X_i^{k-1}g_{ji}X_i^{p-k}-\lambda_{ji}^pg_iX_j,\\
S_{j\ell ,\ell t}=\lambda_{\ell
t}g_{jt}X_{\ell}-\lambda_{ji}X_{\ell}g_{jt}-\lambda_{j\ell}\lambda_{jt}g_{\ell
t}X_j+\lambda_{jt}\lambda_{\ell t}X_tg_{j\ell}.\end{array}$$ This
shows that all overlap elements of $\G$ are reduced to 0 by $\G$,
and hence $\G$ forms a Gr\"obner basis of the ideal
$I=\langle\G\rangle$. Furthermore, it is clear that the requirement
of Corollary 2.5(iii) is satisfied. So the algebra $A=\KS /I$ has a
two-sided Gr\"obner basis theory with respect to the two-sided
admissible systems $(\OV{N(I)},\prec_{_I})$, where $\prec_{_I}$ is
as described in Theorem 2.5(iii). Moreover, in the case where
$\Omega =\left\{\left. g_i=X_i^p~ \right |~1\le i\le n\right\}$, the
algebra $A$ is finite dimensional over $K$. }\par

Note that if $\Omega =\{ g_i=X_i^2~|~1\le i\le n\}$ and all the
$\lambda_{ji}\ne 0$, then the algebra $A$ constructed above is
nothing but the {\it quantum grassmannian} (or {\it quantum
exterior}) algebra in the sense of ([Man], section 3). In this case,
$A$ is clearly $2^n$-dimensional over $K$.{\parindent=0pt\v5

{\bf Remark} Thanks to [G-I1], there are numerous (left, right)
Noetherian algebras of the type described through Theorem 2.2 --
Corollary 2.5, in particular, many {\it quantum binomial algebras}
studied in [Laf] and [G-I2] are such Noetherian algebras, for
instance, the algebras discussed in the foregoing Examples (2) --
(3) with $\lambda ,\mu ,\gamma\ne 0$ and $\gamma^2=1$, respectively
$\lambda^2=1$. By Proposition 1.5, if the algebra $A=R/I$ in Theorem
2.2 -- Corollary 2.5 is respectively left, right, and two-sided
Noetherian, then $A$ has respectively a finite left, right, and
two-sided Gr\"obner basis theory. } \v5

To go further, we let $R=\oplus_{\gamma\in\Gamma}R_{\gamma}$ be a
$\Gamma$-graded $K$-algebra, where $\Gamma$ is an ordered semigroup
by a {\it well-ordering} $<$. For each $f\in R$, say $f=
r_{\gamma_1}+r_{\gamma_2}+\cdots +r_{\gamma_s}$ with
$r_{\gamma_i}\in R_{\gamma_i}$, if $\gamma_1<\gamma_2<\cdots
<\gamma_s$, then we say that $f$ has {\it degree} $\gamma_s$,
denoted $d(f)=\gamma_s$, and we call $r_{\gamma_s}$ the
$\Gamma$-{\it leading homogeneous} element of $f$, denoted $\LH
(f)=r_{\gamma_s}$. Thus, for a subset $S\subset R$ we write $\LH
(S)=\{ \LH (f)~|~f\in S\}$ for the set of $\Gamma$-leading
homogeneous elements of $S$. If $I$ is an ideal of $R$ and $A=R/I$,
then the $\Gamma$-{\it leading homogeneous algebra} of $A$ defined
in [Li2] is the $\Gamma$-graded algebra $A^{\Gamma}_{\rm
LH}(A)=R/\langle\LH (I)\rangle$. Suppose $R$ has a Gr\"obner basis
theory. Our consideration below is to indicate that a given
Gr\"obner basis $\G$ of $I$ may not immediately gives rise to a
Gr\"obner basis theory for $A$ (for instance if it is not the type
considered before), but in case $\LH (\G )$ consists of a subset
$\Omega$ of monomials and $\frac{n(n-1)}{2}$ elements of the form
$g_{ji}$ as in Theorem 2.2, Theorem 2.3, or Theorem 2.4, it does
determine a left, right, or two-sided Gr\"obner basis theory for the
algebra $A^{\Gamma}_{\rm LH}(A)$. The reason that we bring the
object $A^{\Gamma}_{\rm LH}$ into play is sketched as follows (the
interested reader is referred to [Li2] for details).
\par
By ([Li2], Theorem 1.1), there is a $\Gamma$-graded algebra
isomorphism  $A^{\Gamma}_{\rm LH}~\mapright{\cong}{}~G^{\Gamma}(A)$
, where
$G^{\Gamma}(A)=\oplus_{\gamma\in\Gamma}G^{\Gamma}(A)_{\gamma}$ is
the associated $\Gamma$-graded algebra of $A$ defined by the
$\Gamma$-filtration $FA=\{ F_{\gamma}A\}_{\gamma\in\Gamma}$ of $A$
induced by the $\Gamma$-grading filtration $FR=\{
F_{\gamma}R=\oplus_{\gamma'\le \gamma}R_{\gamma
'}\}_{\gamma\in\Gamma}$ of $R$, that is,
$F_{\gamma}A=(F_{\gamma}R+I)/I$,
$G^{\Gamma}(A)_{\gamma}=F_{\gamma}A/F_{\gamma}^*A$ with
$F_{\gamma}^*A=\cup_{\gamma '<\gamma}F_{\gamma '}A$; the
effectiveness of studying $A$ via $A^{\Gamma}_{\rm LH}$ is based on
the diagram
$$\begin{diagram} &&A\\
&\NE&\uTo_{\scriptstyle{\textsf{lifting}}}\\
A^{\Gamma}_{\rm LH}&\rTo^{\cong}&G^{\Gamma}(A)\end{diagram}$$
and the following proposition. {\parindent=0pt\v5

{\bf 2.6. Proposition} ([Li2], Proposition 3.2) Let $R$ and the
notations be as fixed above. Suppose that $R$ has a two-sided
Gr\"obner basis theory with respect to an admissible system $(\B
,\prec_{gr})$, where $\B$ is a {\it skew multiplicative} $K$-{\it
basis} of $R$ consisting of $\Gamma$-{\it homogeneous elements} and
$\prec_{gr}$ is a $\Gamma$-graded two-sided monomial ordering on
$\B$. With notation as before, if $I$ is an ideal of $R$, the
following statements hold.\par (i) Let $\langle\LH (I)\rangle$ be
the $\Gamma$-graded ideal of $R$ generated by the set $\LH (I)$ of
$\Gamma$-leading homogeneous elements of $I$. Then $\LM (I)=\LM
(\langle\LH (I)\rangle )$.\par (ii) A subset $\G\subset I$ is a
Gr\"obner basis of $I$ with respect to $(\B ,\prec_{gr})$ if and
only if $\LH (\G )$, the set of $\Gamma$-leading homogeneous
elements of $\G$, is a Gr\"obner basis for the $\Gamma$-graded ideal
$\langle\LH (I)\rangle$ with respect to $(\B
,\prec_{gr})$.\par\QED}\v5

Proposition 2.6(i) tells us that the set $N(I)$ of normal monomials
in $\B$ (modulo $I$) coincides with the set $N(\langle\LH (I)\rangle
)$ of normal monomials in $\B$ (modulo $\langle\LH (I)\rangle$),
i.e., $N(I)=N(\langle\LH (I)\rangle )$.  Thus we may write
$\widehat{N(I)}=\{ \hat{u}=u+\langle\LH (I)\rangle~|~u\in N(I)\}$
for the $K$-basis of $A^{\Gamma}_{\rm LH}$ determined by
$N(\langle\LH (I)\rangle )$. Note that $\widehat{N(I)}$ consists of
$\Gamma$-homogeneous elements.\par

With notation fixed above, by Proposition 2.6(ii) and Theorem 2.2 --
Theorem 2.4 we are ready to mention the next three results.
{\parindent=0pt\v5

{\bf 2.7. Theorem} Let $R=\oplus_{\gamma\in\Gamma}R_{\gamma}$ and
$(\B ,\prec_{gr})$ be as in Proposition 2.6, and assume further that
$R$ is finitely generated by $\Gamma$-homogeneous elements
$a_1,...,a_n$ (i.e., $R=K[a_1,...,a_n]$), and that $\B$ contains
every monomial of the form
$a_{\ell_1}^{\alpha_1}a_{\ell_2}^{\alpha_2}\cdots
a_{\ell_n}^{\alpha_n}$ with respect to some permutation $a_{\ell_1},
a_{\ell_2}, ..., a_{\ell_n}$ of $a_1, a_2, ..., a_n$, where
$\alpha_1,...,\alpha_n\in\NZ$.  Let $I$ be an ideal of $R$, $A=R/I$,
and $A^{\Gamma}_{\rm LH}=R/\langle\LH (I)\rangle$. If $\G$ is a
minimal Gr\"obner basis of $I$ such that $\LH (\G )=\Omega\cup G$,
where $\Omega\subset\B$ and $G$ consists of $\frac{n(n-1)}{2}$
elements
$$\begin{array}{rcl} g_{ji}&=&a_{\ell_j}a_{\ell_i}-\lambda_{ji}a_{\ell_i}a_{\ell_p},~
1\le i<j\le n,~ i<p\le n,~\lambda_{ji}\in K^*\cup\{ 0\},\\
&{~}&\hbox{satisfying}~\LM (g_{ji})=a_{\ell_j}a_{\ell_i},~1\le
i<j\le n,\end{array}$$ then the following statements hold.\par

(i) $N(I)\subseteq \{
a_{\ell_1}^{\alpha_1}a_{\ell_2}^{\alpha_2}\cdots
a_{\ell_n}^{\alpha_n}~|~\alpha_1,...,\alpha_n\in\NZ\}$, and
$\widehat{N(I)}$ is a skew multiplicative $K$-basis for
$A^{\Gamma}_{\rm LH}$.\par

(ii)  Let $\prec_{_{\LH (I)}}$ denote the ordering on
$\widehat{N(I)}$, which is induced by the lexicographic ordering or
the $\NZ$-graded lexicographic ordering on $\NZ^n$, such that
$$\widehat{a_{\ell_n}}\prec_{_{\LH (I)}}
\widehat{a_{\ell_{n-1}}}\prec_{_{\LH (I)}}\cdots \prec_{_{\LH
(I)}}\widehat{a_{\ell_2}}\prec_{_{\LH (I)}}\widehat{a_{\ell_1}}.$$
If
$\widehat{a_{\ell_1}}^{\alpha_1}\widehat{a_{\ell_2}}^{\alpha_2}\cdots
\widehat{a_{\ell_n}}^{\alpha_n}\ne 0$ implies
$a_{\ell_1}^{\alpha_1}a_{\ell_2}^{\alpha_2}\cdots
a_{\ell_n}^{\alpha_n}\in N(I)$, then $\prec_{_{\LH (I)}}$ is a left
monomial ordering on $\widehat{N(I)}$, and hence $A^{\Gamma}_{\rm
LH}$ has a left Gr\"obner basis theory with respect to the left
admissible system $(\widehat{N(I)},\prec_{_{\LH (I)}} )$.
\par\QED\v5

{\bf 2.8. Theorem} Under the assumptions as in Theorem 2.7, let $I$
be an ideal of $R$, $A=R/I$, and $A^{\Gamma}_{\rm LH}=R/\langle\LH
(I)\rangle$. If $\G$ is a minimal Gr\"obner basis of $I$ such that
$\LH (\G )=\Omega\cup G$, where $\Omega\subset\B$ and $G$ consists
of $\frac{n(n-1)}{2}$ elements
$$\begin{array}{rcl} g_{ji}&=&a_{\ell_j}a_{\ell_i}-\lambda_{ji}a_{\ell_p}a_{\ell_j},~
1\le i<j\le n,~ p<j,~\lambda_{ji}\in K^*\cup\{ 0\},\\
&{~}&\hbox{satisfying}~\LM (g_{ji})=a_{\ell_j}a_{\ell_i},~1\le
i<j\le n,\end{array}$$ then the following statements hold.\par

(i) $N(I)\subseteq \{
a_{\ell_1}^{\alpha_1}a_{\ell_2}^{\alpha_2}\cdots
a_{\ell_n}^{\alpha_n}~|~\alpha_1,...,\alpha_n\in\NZ\}$, and
$\widehat{N(I)}$ is a skew multiplicative $K$-basis for
$A^{\Gamma}_{\rm LH}$.\par

(ii)  Let $\prec_{_{\LH (I)}}$ denote the ordering on
$\widehat{N(I)}$, which is induced by the $\NZ$-graded reverse
lexicographic ordering on $\NZ^n$, such that
$$\widehat{a_{\ell_n}}\prec_{_{\LH (I)}}
\widehat{a_{\ell_{n-1}}}\prec_{_{\LH (I)}}\cdots \prec_{_{\LH
(I)}}\widehat{a_{\ell_2}}\prec_{_{\LH (I)}}\widehat{a_{\ell_1}}.$$
If
$\widehat{a_{\ell_1}}^{\alpha_1}\widehat{a_{\ell_2}}^{\alpha_2}\cdots
\widehat{a_{\ell_n}}^{\alpha_n}\ne 0$ implies
$a_{\ell_1}^{\alpha_1}a_{\ell_2}^{\alpha_2}\cdots
a_{\ell_n}^{\alpha_n}\in N(I)$, then $\prec_{_{\LH (I)}}$ is a right
monomial ordering on $\widehat{N(I)}$, and hence $A^{\Gamma}_{\rm
LH}$ has a right Gr\"obner basis theory with respect to the right
admissible system $(\widehat{N(I)},\prec_{_{\LH (I)}} )$.
\par\QED\v5

{\bf 2.9. Theorem} Under the assumptions as in Theorem 2.7, let $I$
be an ideal of $R$, $A=R/I$, and $A^{\Gamma}_{\rm LH}=R/\langle\LH
(I)\rangle$. If $\G$ is a minimal Gr\"obner basis of $I$ such that
$\LH (\G )=\Omega\cup G$, where $\Omega\subset\B$ and $G$ consists
of $\frac{n(n-1)}{2}$ elements
$$\begin{array}{rcl} g_{ji}&=&a_{\ell_j}a_{\ell_i}-\lambda_{ji}a_{\ell_i}a_{\ell_j},~
1\le i<j\le n,~ \lambda_{ji}\in K^*\cup\{ 0\},\\
&{~}&\hbox{satisfying}~\LM (g_{ji})=a_{\ell_j}a_{\ell_i},~1\le
i<j\le n,\end{array}$$                            Then the following
statements hold.\par

(i) $N(I)\subseteq \{
a_{\ell_1}^{\alpha_1}a_{\ell_2}^{\alpha_2}\cdots
a_{\ell_n}^{\alpha_n}~|~\alpha_1,...,\alpha_n\in\NZ\}$, and
$\widehat{N(I)}$ is a skew multiplicative $K$-basis for
$A^{\Gamma}_{\rm LH}$\par

(ii) If
$\widehat{a_{\ell_1}}^{\alpha_1}\widehat{a_{\ell_2}}^{\alpha_2}\cdots
\widehat{a_{\ell_n}}^{\alpha_n}\ne 0$ implies
$a_{\ell_1}^{\alpha_1}a_{\ell_2}^{\alpha_2}\cdots
a_{\ell_n}^{\alpha_n}\in N(I)$, then $A^{\Gamma}_{\rm LH}$ has a
two-sided Gr\"obner basis theory with respect to the admissible
system $(\widehat{N(I)},\prec_{_{\LH (I)}} )$, where $\prec_{_{\LH
(I)}} $ is the two-sided monomial ordering on $\widehat{N(I)}$
induced by  one of the following orderings on $\NZ^n$: the
lexicographic ordering, the $\NZ$-graded lexicographic ordering, and
the $\NZ$-graded reverse lexicographic ordering,  such that
$$\widehat{a_{\ell_n}}\prec_{_{\LH (I)}} \widehat{a_{\ell_{n-1}}}\prec_{_{\LH (I)}} \cdots
\prec_{_{\LH (I)}} \widehat{a_{\ell_2}}\prec_{_{\LH (I)}}
\widehat{a_{\ell_1}}.$$
\par\QED}\v5

Examples given below not only illustrate Theorem 2.7 -- 2.9, but
also answers why the three theorems said nothing about the existence
of a (left, right, or two-sided) Gr\"obner basis theory for the
quotient algebra $A=R/I$. {\parindent=0pt \v5

{\bf Example} (5) Let $\G$ be the subset of the free $K$-algebra
$\KS =K\langle X_1,X_2,X_3\rangle$ consisting of
$$\begin{array}{l} g_{21}=X_2X_1-\lambda X_1X_3+\alpha X_1,\\
g_{31}=X_3X_1,\\
g_{32}=X_3X_2-\mu X_2X_3+\alpha X_3,\end{array}~\hbox{where}~
\lambda ,\mu ,\alpha\in K^*\cup\{ 0\}.$$ Then, with respect to the
natural $\NZ$-graded lexicographic ordering subject to
$X_1\prec_{grlex}X_2\prec_{grlex}X_3$, where each  $X_i$ is assigned
to the degree 1, it is straightforward to verify that $\LM
(g_{ji})=X_jX_i$, $1\le i<j\le 3$, $\G$ forms a Gr\"obner basis of
the ideal $I=\langle\G\rangle$ and $N(I)=\{
X_1^{\alpha_1}X_2^{\alpha_2}X_3^{\alpha_3}~|~\alpha_j\in\NZ\}$.
Consider the quotient algebra $A=\KS /I$. Since $\G$ is of the type
required by Theorem 2.7 such that $\LH (\G)=\{X_2X_1-\lambda
X_1X_3,~X_3X_1,~X_3X_2-\mu X_2X_3\}$, it follows that the algebra
$A^{\NZ}_{\rm LH}=\KS /\langle\LH (\G )\rangle$ has a left Gr\"obner
basis theory with respect to the left admissible system
$(\widehat{N(I)},\prec_{_{\LH (I)}})$, where $\prec_{_{\LH (I)}}$ is
as described in Theorem 2.7.}\par

If  $\mu =0$ and $\alpha\ne 0$, then note that the algebra $A$ has
$\OV{X_2}$, $\OV{X_2}^2$, $\OV{X_3}\in\OV{N(I)}$ such that
$\OV{X_2}\ne\OV{X_2}^2$, $\OV{X_3}\cdot\OV{X_2}=-\alpha\OV{X_3}$,
and $\OV{X_3}\cdot\OV{X_2}^2=-\alpha\OV{X_3}\cdot\OV{X_2}$. Thus,
with respect to any total ordering $\prec$ on $\OV{N(I)}$, we have
$$\LM (\OV{X_3}\cdot\OV{X_2}^2)=\OV{X_3}=\LM
(\OV{X_3}\cdot\OV{X_2}).$$ This shows that $\prec$ cannot be a
(left, right, or two-sided) monomial ordering on $\OV{N(I)}$, and
therefore, the algebra $A$ cannot have a (left, right, and
two-sided) Gr\"obner basis whenever the $K$-basis $\OV{N(I)}$ is
used.{\parindent=0pt\v5

{\bf Example} (6) Let $\G$ be the subset of the free $K$-algebra
$\KS =K\langle X_1,X_2,X_3\rangle$ consisting of
$$\begin{array}{l} g_{21}=X_2X_1-\lambda X_1X_2+\alpha X_1,\\
g_{31}=X_3X_1,\\
g_{32}=X_3X_2-\mu X_1X_3+\alpha X_3,\end{array}~\hbox{where}~
\lambda ,\mu ,\alpha\in K^*\cup\{ 0\}.$$ Then, with respect to the
natural $\NZ$-graded lexicographic ordering subject to
$X_1\prec_{grlex}X_2\prec_{grlex}X_3$, where each  $X_i$ is assigned
to the degree 1, it is straightforward to verify that $\LM
(g_{ji})=X_jX_i$, $1\le i<j\le 3$, $\G$ forms a Gr\"obner basis of
the ideal $I=\langle\G\rangle$ and $N(I)=\{
X_1^{\alpha_1}X_2^{\alpha_2}X_3^{\alpha_3}~|~\alpha_j\in\NZ\}$.
Consider the quotient algebra $A=\KS /I$. Since $\G$ is of the type
required by Theorem 2.8 such that $\LH (\G)=\{X_2X_1-\lambda
X_1X_2,~X_3X_1,~X_3X_2-\mu X_1X_3\}$, it follows that the algebra
$A^{\NZ}_{\rm LH}=\KS /\langle\LH (\G )\rangle$ has a right
Gr\"obner basis theory with respect to the right admissible system
$(\widehat{N(I)},\prec_{_{\LH (I)}})$, where $\prec_{_{\LH (I)}}$ is
as described in Theorem 2.8.}\par

If  $\lambda =0$ and $\alpha\ne 0$, then note that the algebra $A$
has $\OV{X_2}$, $\OV{X_2}^2$, $\OV{X_1}\in\OV{N(I)}$ such that
$\OV{X_2}\ne\OV{X_2}^2$, $\OV{X_2}\cdot\OV{X_1}=-\alpha\OV{X_1}$,
and $\OV{X_2}^2\cdot\OV{X_1}=-\alpha\OV{X_2}\cdot\OV{X_1}$. Thus,
with respect to any total ordering $\prec$ on $\OV{N(I)}$, we have
$$\LM (\OV{X_2}^2\cdot\OV{X_1})=\OV{X_1}=\LM
(\OV{X_2}\cdot\OV{X_1}).$$ This shows that $\prec$ cannot be a
(left, right, or two-sided) monomial ordering on $\OV{N(I)}$, and
therefore, the algebra $A$ cannot have a (left, right, or two-sided)
Gr\"obner basis whenever the $K$-basis $\OV{N(I)}$ is
used.{\parindent=0pt\v5

{\bf Example} (7) Let $\G$ be the subset of the free $K$-algebra
$\KS =K\langle X_1,X_2,X_3\rangle$ consisting of
$$\begin{array}{l} g_{31}=X_3X_1-\lambda X_1X_3+\gamma X_3,\\
g_{12}=X_1X_2-\lambda X_2X_1+\gamma X_2,\\
g_{32}=X_3X_2-\omega X_2X_3+f(X_1),\end{array}$$ where $\lambda
,\gamma ,\omega\in K^*\cup\{ 0\}$, and $f(X_1)$ is a polynomial in
the variable $X_1$. If $d(f)\le 2$, we assign
$d(X_1)=d(X_2)=d(X_3)=1$; if $d(f)=n\ge 3$, we assign
$d(X_1)=d(X_2)=1$, $d(X_3)=n$. Then, using the $\NZ$-graded
lexicographic ordering $X_2\prec_{grlex}X_1\prec_{grlex}X_3$ with
respect to either  $\NZ$-gradation assigned to $\KS$ above, it is
straightforward to verify that $\LM (g_{ji})=X_jX_i$, $1\le i<j\le
3$, $\G$ forms a Gr\"obner basis of the ideal $I=\langle\G\rangle$
and $N(I)=\{
X_2^{\alpha_2}X_1^{\alpha_1}X_3^{\alpha_3}~|~\alpha_j\in\NZ\}$. By
Proposition 2.6(ii), $$\begin{array}{l} \LH (\G )=\{X_3X_1-\lambda
X_1X_3,~X_1X_2-\lambda X_2X_1,~X_3X_2-\omega
X_2X_3\}~\hbox{if}~d(f)\ge 3,\\ \hbox{respectively}\\
\LH (\G )=\{X_3X_1-\lambda X_1X_3,~X_1X_2-\lambda
X_2X_1,~X_3X_2-\omega X_2X_3+aX_1^2\} \\
~~~~~~~~~~~~~~\hbox{if}~f(x_1)=aX_1^2+bX_1+c~\hbox{with}~a,b,c\in K,
\end{array}$$ is a Gr\"obner basis for the ideal  $\langle\LH (I)\rangle$ with
respect to the respective $\NZ$-gradation assigned to $\KS$. The
existence of a Gr\"obner basis theory for the algebras $A=\KS /I$
and $A^{\NZ}_{\rm LH}=\KS/\langle\LH (\G )\rangle$ is discussed as
follows.
\par

(1) If $\lambda\ne 0$ and $\omega\ne 0$, then, in any case (i.e.,
$d(f)\le 2$ or $d(f)\ge 3)$,  with respect to
$\OV{X_2}\prec_{grlex}\OV{X_1}\prec_{grlex}\OV{X_3}$ on $\OV{N(I)}$,
respectively
$\widehat{X_2}\prec_{grlex}\widehat{X_1}\prec_{grlex}\widehat{X_3}$
on $\widehat{N(I)}$, both $A$ and $A^{\NZ}_{\rm LH}$ are solvable
polynomial algebras in the sense of [K-RW]. Hence,  every (left,
right, or two-sided) ideal has a finite  Gr\"obner basis.\par

(2) If $d(f)\ge 3$ and $\lambda =0$ or $\omega =0$, then by Theorem
2.9, $A^{\NZ}_{\rm LH}$ has a two-sided Gr\"obner basis theory with
respect to the admissible system $(\widehat{N(I)},\prec_{_{\LH
(I)}})$, where $\prec_{_{\LH (I)}}$ is the monomial ordering as
described in Theorem 2.9.\par

(3) In the case that $f(X_1)=aX_1^2+bX_1+c=0$, if $a=0$, then, for
arbitrary $\lambda$ and $\omega$,  $A^{\NZ}_{\rm LH}$ has a
two-sided Gr\"obner basis theory with respect to the admissible
system $(\widehat{N(I)},\prec_{_{\LH (I)}})$, where $\prec_{_{\LH
(I)}}$ is the monomial ordering as described in Theorem 2.9.\par

(4) In the case that $f(X_1)=aX_1^2+bX_1+c=0$, if $a\ne 0$ and
$\lambda =0$ or $\omega =0$, then we are not clear about the
existence of a (left, right, or two-sided) gr\"obner basis for
$A^{\NZ}_{\rm LH}$.\par

(5) If  $\lambda =0$ and $\gamma\ne 0$, then, in any case (i.e.,
$d(f)\le 2$ or $d(f)\ge 3)$, note that $\OV{X_1}$, $\OV{X_1}^2$,
$\OV{X_3}\in\OV{N(I)}$, $\OV{X_1}\ne\OV{X_1}^2$,
$\OV{X_3}\cdot\OV{X_1}=-\gamma\OV{X_3}$, and
$\OV{X_3}\cdot\OV{X_1}^2=\gamma^2\OV{X_3}$. Thus,  with respect to
any total ordering $\prec$ on $\OV{N(I)}$, we have
$$\LM (\OV{X_3}\cdot\OV{X_1})=\OV{X_3}=\LM
(\OV{X_3}\cdot\OV{X_1}^2).$$ This shows that $\prec$ cannot be a
(left, right, or two-sided) monomial ordering on $\OV{N(I)}$, and
therefore, the quotient algebra $A=R/I$ cannot have a (left, right,
or two-sided) Gr\"obner basis theory whenever the $K$-basis
$\OV{N(I)}$ is used.\v5

{\bf Remark} Here let us point out that Example (7) discussed above
applies to two interesting cases.
\par

(1) In the case where $f(X_1)=aX_1^2+bX_1$ such that
$\lambda\gamma\omega b\ne 0$, the algebra $A=R/I$ is nothing but the
{\it conformal {\rm sl$_2$} enveloping algebra} in the sense of
[Le-B]. \par

(2) Suppose that the field $K$ is algebraically closed and
$f(X_1)=X_1$. If  $\lambda$, $\omega$ are roots of the equation
$x^2-\alpha x-\beta =0$, then since $\lambda +\omega =\alpha$,
$-\lambda\omega =\beta$, it follows from [KMP], [CM] and [CS] that
all {\it down-up algebras} $A(\alpha ,\beta ,\gamma )$ in the sense
of [Ben], are recaptured by the quadric defining relations
$$\begin{array}{l} g_{31}=X_3X_1-\lambda X_1X_3+\gamma X_3,\\
g_{12}=X_1X_2-\lambda X_2X_1+\gamma X_2,\\
g_{32}=X_3X_2-\omega X_2X_3+X_1.\end{array}$$ }\v5

In addition to the effectiveness of studying $A$ via a Gr\"obner
basis theory of $A^{\Gamma}_{\rm LH}$ (as may be seen from [Li2]),
the natural question is now to ask if there would be any appropriate
condition on $A$ or on $A^{\Gamma}_{\rm LH}$ that enables us to
obtain a Gr\"obner basis theory for $A$ subject to a Gr\"obner basis
theory of $A^{\Gamma}_{\rm LH}$.\par Enlightened by ([Li2], Theorem
1.1), the next theorem, which may be viewed as a recognizable
version of ([Li1], CH.IV, Theorem 1.6), provides a solution to the
question posed above. Before mentioning the next theorem, we need a
little more preparation, for, the ideal $I$ concerned by the theorem
is not necessarily generated by monomials and binomials, and
moreover, the $K$-basis $\OV{N(I)}$ of $A=R/I$, respectively the
$K$-basis $\widehat{N(I)}$ of $A^{\Gamma}_{\rm LH}$, is not
necessarily a skew multiplicative $K$-basis. In this case, with any
fixed $K$-basis $\mathscr{B}$ of $A$, the definition of a (left,
right, or two-sided) monomial ordering $\prec$ on $\mathscr{B}$ is
the same as that defined in Definition 1.1; for $u,v\in\mathscr{B}$,
we say that $u$ divides $v$, denoted $u|v$, if there are
$w,s\in\mathscr{B}$ such that $v=\LM (wus)$; and we say that $A$ has
a (left, right, or two-sided) Gr\"obner basis theory (comparing with
Definition 1.4) if, with respect to the (left, right, or two-sided)
admissible system  $(\mathscr{B},\prec )$, every proper (left,
right, or two-sided) ideal $I$ of $A$ has a generating set $\G$ such
that $f\in I$ implies $\LM (g)|\LM (f)$ for some $g\in\G$. Similar
convention is valid for $A^{\Gamma}_{\rm LH}$. {\parindent=0pt\v5

{\bf 2.10. Theorem} Let $R=\oplus_{\gamma\in\Gamma}R_{\gamma}$ and
$(\B ,\prec_{gr})$ be as in Proposition 2.6, and let $I$ be an
arbitrary proper ideal of $R$ (i.e., $I$ is not necessarily
generated by monomials and binomials). Put $A=R/I$, $A^{\Gamma}_{\rm
LH}=R/\langle\LH (I)\rangle$. With notation as before and the
convention made above, suppose that
\parindent=.7truecm \par
\item{(a)} $A^{\Gamma}_{\rm LH}$ is a domain; and \par
\item{(b)} $A^{\Gamma}_{\rm LH}$ has a (left, right, or two-sided) Gr\"obner basis theory
with respect to some  (left, right, or two-sided) admissible system
$(\widehat{N(I)},\prec_{\widehat{gr}})$, where
$\prec_{\widehat{gr}}$ is a $\Gamma$-graded (left, right, or
two-sided) monomial ordering  on $\widehat{N(I)}$, such that every
$\Gamma$-graded (left, right, or two-sided) ideal of
$A^{\Gamma}_{\rm LH}$ has a (left, right, or two-sided) Gr\"obner
basis consisting of $\Gamma$-homogeneous elements.
\parindent=0pt\par  Then the order $\prec$ on $\OV{N(I)}$ defined
subject to the rule: for $\OV u$, $\OV v\in\OV{N(I)}$,
$$\bar{u}\prec\bar{v}~\hbox{whenever}~\hat{u}\prec_{\widehat{gr}}\hat{v},$$
is a (left, right, or two-sided) monomial ordering on $\OV{N(I)}$,
such that every (left, right, or two-sided) ideal of $A$ has a
(left, right, or two-sided) Gr\"obner basis. \vskip 6pt

{\bf Proof} Note that it is not difficult to check that if $R$ has a
left or right Gr\"obner basis theory with respect to some left or
right $\Gamma$-graded monomial ordering on $\B$, then an analogue of
Proposition 2.6 holds for a left or right ideal. The proof given
below will show that it  is enough to conclude that $A$ has a
two-sided Gr\"obner basis theory whenever $A^{\Gamma}_{\rm LH}$ has
a two-sided Gr\"obner basis theory.}
\par

Since $\prec$ is clearly a well-ordering on $\OV{N(I)}$, we show
first that $\prec$ satisfies the axioms (MO1) -- (MO2) of Definition
1.1(i). We start by finding the relation between $\LM (\hat{f}~)$
and $\LM (\bar{f}~)$ for a nonzero $f\in R$, where
$\hat{f}=f+\langle\LH (I)\rangle\in A^{\Gamma}_{\rm LH}$,
$\bar{f}=f+I\in A$. First note that $I^*=I-\{ 0\}$ is trivially a
Gr\"obner basis of $I$, and hence $\LH (I^*)$ is a Gr\"obner basis
for $\langle\LH (I)\rangle$ by Proposition 2.6(ii). Thus, writing
$f=\LH (f)+f'$ with $d(f')<d(f)=d(\LH (f))$ and noticing that
$\prec_{gr}$ is a $\Gamma$-graded monomial ordering on $\B$, the
division by $\LH (I^*)$ turns out that the homogeneous element $\LH
(f)$ has a Gr\"obner presentation (Definition 1.2):
$$\begin{array}{rcl} \LH (f)&=&\SUM^{}_{i,j}\lambda_{ij}w_{ij}\LH (f_j)v_{ij}+
\SUM^{}_{\ell}\lambda_{\ell}u_{\ell}~\hbox{with}~\lambda_{ij}\in
K^*,~\lambda_{\ell}\in
K,~w_{ij},v_{ij}\in\B,~f_j\in I^*, \\
&{~}&\hbox{where if}~\lambda_{\ell}\ne 0~\hbox{then}~~u_{\ell}\in
N(\langle\LH (I)\rangle )=N(I),~d(u_{\ell})=d(\LH
(f))=d(f),\end{array}$$
and that if $f'\ne 0$ then $f'$ has a Gr\"obner presentation
$$\begin{array}{rcl} f'&=&\SUM^{}_{p,q}\lambda_{pq}w_{pq}\LH (f_q)v_{pq}+
\SUM^{}_{h}\lambda_{h}u_{h}~\hbox{with}~\lambda_{pq}\in
K^*,\lambda_{h}\in
K,~w_{pq},v_{pq}\in\B,~f_q\in I^*, \\
&{~}&\hbox{where if}~\lambda_h\ne 0~\hbox{then}~u_{h}\in
N(\langle\LH (I)\rangle )=N(I),~d(u_{h})<d(\LH
(f))=d(f).\end{array}$$
Consequently,
$$\begin{array}{rcl} f=\LH (f)+f'&=&\SUM^{}_{i,j}\lambda_{ij}w_{ij}\LH (f_j)v_{ij}+
\SUM^{}_{p,q}\lambda_{pq}w_{pq}\LH
(f_q)v_{pq}\\
\\
&{~}&~+\SUM^{}_{d(u_{\ell})=d(f)}\lambda_{\ell}u_{\ell}+
\SUM^{}_{d(u_h)<d(f)}\lambda_{h}u_{h}.\end{array}$$
If we write each $f_j$ appearing in the presentation of $\LH (f)$
obtained above as $f_j=\LH (f_j)+f_j'$ with $d(f_j')<d(f_j)$, then,
since $\Gamma$ is an ordered semigroup, $f$ has a presentation
$$\begin{array}{rcl} f&=&\SUM^{}_{i,j}\lambda_{ij}w_{ij}f_jv_{ij}+
f'-\SUM^{}_{ij}\lambda_{ij}w_{ij}f_j'v_{ij}+\SUM^{}_{\ell}\lambda_{\ell}u_{\ell}~~
\hbox{satisfying}\\
\\
&{~}&d\left (f'-\SUM^{}_{ij}\lambda_{ij}w_{ij}f_j'v_{ij}\right )
<d(f),~u_{\ell}\in N(\langle\LH (I)\rangle
)=N(I)~\hbox{and}~d(u_{\ell})=d(f).\end{array}$$
If furthermore we do division by $I^*$ with respect to $\prec_{gr}$,
then the remainder of $f'-\SUM^{}_{ij}\lambda_{ij}w_{ij}f_j'v_{ij}$
must have degree $< d(f)$. Summing up, recalling that the
$\Gamma$-gradation of $A^{\Gamma}_{\rm LH}$ is the one induced by
the $\Gamma$-gradation of $R$, and that the order $\prec$ on
$\OV{N(I)}$ is defined subject to the $\Gamma$-graded monomial
ordering $\prec_{\widehat{gr}}$ on $\widehat{N(I)}$, we are now able
to conclude that if $\LH (f)\not\in\langle\LH (I)\rangle$, or
equivalently, if $\widehat{\LH (f)}\ne 0$, then
$$\LM (\hat{f}~)=\hat{u}~\hbox{implies}~\LM (\bar{f}~)=\bar{u},~\hbox{where}~u\in N(I).\eqno{(1)}$$ \par
Next, for $\bar{u}$, $\bar{v}$, $\bar{w}$, $\bar{s}\in\OV{N(I)}$,
suppose $\bar{u}\prec\bar{v}$, $\LM (\bar{w}\bar{u}\bar{s})\not\in
K^*\cup\{ 0\}$ and $\LM (\bar{w}\bar{v}\bar{s})\not\in K^*\cup\{
0\}$. Then $\hat{u}\prec_{\widehat{gr}}\hat{v}$ by the definition of
$\prec$. Also we conclude that $\LM (\hat{w}\hat{u}\hat{s})\not\in
K^*\cup\{ 0\}$, $\LM (\hat{w}\hat{v}\hat{s})\not\in K^*\cup\{ 0\}$,
and hence
$$\LM (\widehat{wus})=\LM (\hat{w}\hat{u}\hat{s})\prec_{\widehat{gr}}
\LM (\hat{w}\hat{v}\hat{s})=\LM (\widehat{wvs}).\eqno{(2)}$$
To see this, put $f=wus$, $g=wvs$. Then, noticing that
$w,u,v,s\in\B$ are $\Gamma$-homogeneous elements of $R$, we have
$f=\LH (f)=wus$, $g=\LH (g)=wvs$, and thereby $\hat{f}=\widehat{\LH
(f)}=\widehat{wus}\ne 0$, $\hat{g}=\widehat{\LH
(g)}=\widehat{wvs}\ne 0$, since $A^{\Gamma}_{\rm LH}$ is a domain by
the assumption (a).  Let $\LM (\hat{f}~)=\widehat{u_1}$ and $\LM
(\hat{g}~)=\widehat{v_1}$ with $u_1$, $v_1\in N(I)$. Then it follows
from the assertion (1) above that
$$\begin{array}{l} \OV{u_1}=\LM (\bar{f}~)=\LM (\OV{wus})=\LM (\bar{w}\bar{u}\bar{s}),\\
\OV{v_1}=\LM (\bar{g}~)=\LM (\OV{wvs})=\LM
(\bar{w}\bar{v}\bar{s}).\end{array}\eqno{(3)}$$ Thus $\LM
(\hat{w}\hat{u}\hat{s})\not\in K^*\cup\{ 0\}$, $\LM
(\hat{w}\hat{v}\hat{s})\not\in K^*\cup\{ 0\}$ and (2) follows.
Consequently, (2) $+$ (3) gives rise to
$$\LM (\bar{w}\bar{u}\bar{s})=\OV{u_1}\prec\OV{v_1}=\LM (\bar{w}\bar{v}\bar{s}).$$
This proves that $\prec$ satisfies (MO1). It remains to show that
$\prec$ satisfies (MO2) as well. Suppose that $\bar{u}$, $\bar{w}$,
$\bar{v}$, $\bar{s}\in\OV{N(I)}$ satisfy $\bar{u}=\LM
(\bar{w}\bar{v}\bar{s})$ (where $\bar{w}\ne\bar{1}$ or
$\bar{s}\ne\bar{1}$ if $1\in\B$, hence $w\ne 1$ or $s\ne 1$ in
$N(I)$ and  $\hat{w}\ne \hat{1}$ or $\hat{s}\ne \hat{1}$ in
$\widehat{N(I)}$). Since $A^{\Gamma}_{\rm LH}$ is a domain, $\LM
(\hat{w}\hat{v}\hat{s})\ne 0$. So, if $\LM
(\hat{w}\hat{v}\hat{s})=\LM (\widehat{wvs})=\widehat{v_1}$ with
$v_1\in N(I)$, then, as argued above, it follows from the foregoing
assertion (1) that $\bar{u}=\LM (\bar{w}\bar{v}\bar{s})=\LM
(\OV{wvs})=\OV{v_1}$. Thus,
$\hat{v}\prec_{\widehat{gr}}\widehat{v_1}$ implies that
$\bar{v}\prec\bar{u}$, as desired.\par

Finally, we  prove that if $\OV J=J/I$  is a proper ideal of
$A=R/I$, where $J$ is an ideal of $R$ containing $I$, then $\OV J$
has a Gr\"obner basis $\OV{\G}$ with respect to the admissible
system $(\prec ,\OV{N(I)})$ obtained above. To this end, let us note
first that if $f\in R-I$, then since $\Gamma$ is ordered by a
well-ordering and $\langle\LH (I)\rangle$ is a $\Gamma$-graded ideal
of $R$, there is some $f'\in R$ such that $\LH (f')\not\in\langle\LH
(I)\rangle$ and $\bar{f}=\OV{f'}$. So, without loss of generality,
we may always assume that
$$f\in R-I~\hbox{implies}~\LH (f)\not\in\langle\LH (I)\rangle .\eqno{(4)}$$
Consequently, since the $\Gamma$-gradation of $A^{\Gamma}_{\rm LH}$
is induced by the $\Gamma$-gradation of $R$, and since
$\prec_{\widehat{gr}}$ is a $\Gamma$-graded monomial ordering on
$\widehat{N(I)}$, if we write $f=\LH (f)+f_1$ then
$$\begin{array}{l} \LH (\hat{f})=\widehat{\LH (f)},\\
\LM (\hat{f})=\LM (\LH (\hat{f}))=\LM \left (\widehat{\LH
(f)}\right).\end{array}\eqno{(5)}$$ Now, noticing that $\LH
(I)\subset \LH (J)$ and thus $\widehat{J}=\langle\LH (J)\rangle
/\langle\LH(I)\rangle$ is a $\Gamma$-graded ideal of
$A^{\Gamma}_{\rm LH}$, by Proposition 1.5 and Proposition 2.6(i),
$\widehat{J}\ne A^{\Gamma}_{\rm LH}$. Hence, with respect to the
data $(\widehat{N(I)},\prec_{\widehat{gr}})$, $\widehat{J}$ has a
Gr\"obner basis $\mathscr{G}$ consisting of $\Gamma$-homogeneous
elements. Once again since the $\Gamma$-gradation of
$A^{\Gamma}_{\rm LH}$ is induced by the $\Gamma$-gradation of $R$,
and thus $\widehat{\LH (J)}=\{\widehat{\LH (h)}~ |~h\in J\}$ is a
generating set of $\widehat{J}$ consisting of $\Gamma$-homogeneous
elements, it is straightforward to check (or see Lemma 4.2 in later
Section 4) that $\mathscr{G}\subset\widehat{\LH (J)}$. Hence there
is a subset $\G\subset J-I$ such that $\widehat{\LH (\G
)}=\mathscr{G}$. Put $\OV{\G}=\{\bar{g}~|~g\in\G\}$. If
$\bar{f}\in\OV{J}-\{ 0\}$, then by the previous (4) and (5), there
is some $g\in\G$ such that $\LM (\widehat{\LH (g)})|\LM (\bar{f})$,
and consequently $\LM (\hat{g})|\LM (\hat{f})$, i.e., there are
$\hat{u},\hat{v}\in\widehat{N(I)}$ such that $\LM (\hat{f})=\LM
(\hat{u}\LM (\hat{g})\hat{v})$. It follows from the foregoing (1)
and the argument given before the formula (3) that $$\LM
(\bar{f})=\LM (\bar{u}\LM (\bar{g})\bar{v}).$$ This shows that
$\OV{\G}$ is a Gr\"obner basis for $\OV{J}$. \QED{\parindent=0pt\v5

{\bf Remark} From the proof given above it is clear that if every
$\Gamma$-graded (left, right, or two-sided) ideal of
$A^{\Gamma}_{\rm LH}$ has a finite (left, right, or two-sided)
Gr\"obner basis consisting of $\Gamma$-homogeneous elements. Then
every (left, right, or two-sided) ideal of $A$ has a finite (left,
right, or two-sided) Gr\"obner basis.}\v5

There are examples of algebras other than algebras of the solvable
type in the sense of [K-RW] to illustrate Theorem 2.10.
{\parindent=0pt\v5

{\bf Example} (8) If in Corollary 2.5(i) -- (ii) we have $\Omega
=\emptyset$ and all the $\lambda_{ji}\ne 0$, then it follows from
the proof of Theorem 2.2 (also a similar proof of Theorem 2.3,
though it is not mentioned there) that the quadratic algebra
considered in (i) -- (ii) respectively is a domain. Hence, by
Proposition 2.6, if $\G$ is a Gr\"obner basis of the ideal
$I=\langle\G\rangle$ in $K\langle X_1,...,X_n\rangle$ such that $\LH
(\G )$ is of the form as described in  Corollary 2.5(i) -- (ii)
respectively but with $\Omega =\emptyset$ and all the
$\lambda_{ji}\ne 0$, then the quotient algebra $K\langle
X_1,...,X_n\rangle /I$ will provide us with a practical algebra
required by Theorem 2.10.  Omitting tedious verification, we list a
family of the desired Gr\"obner bases $\G$ in $K\langle
X_1,X_2,X_3,X_4\rangle$ consisting of
$$\begin{array}{ll} g_{43}=X_4X_3-\lambda_{43}X_3X_4-\alpha
X_2,&g_{32}=X_3X_2-
\lambda_{32}X_2X_3-\gamma X_4,\\
g_{42}=X_4X_2- \lambda_{42}X_2X_4+\beta X_3,&g_{31}=X_3X_1-
\lambda_{31}X_1X_4,\\
g_{41}=X_4X_1- \lambda_{41}X_1X_2,&g_{21}=X_2X_1-
\lambda_{21}X_1X_3,\end{array}$$
where $\lambda_{ji},\alpha ,\beta ,\gamma\in K^*$ satisfying
$$\begin{array}{l} \lambda_{21}\alpha =\beta,~
\lambda_{31}\beta =\gamma ,~\lambda_{41}\gamma =\alpha ,\\
\lambda_{43}\lambda_{42}=1,~\lambda_{43}=\lambda_{32},
\end{array}$$
and the monomial ordering used here is the $\NZ$-graded
lexicographic monomial ordering $\prec_{grlex}$ with respect to
$d(X_i)=1$, $1\le i\le 4$, such that
$X_1\prec_{grlex}X_2\prec_{grlex}X_3\prec_{grlex}X_4$.}

\section*{3. Skew 2-nomial Algebras}
In Section 2 we have listed some practical algebras (including
numerous quantum binomial algebras in the sense of [Laf] and [G-I2])
defined by elements of the form $u-\lambda v$ and $w$, where $u,v,w$
are monomials and $\lambda\in K^*$, which are not the familiar types
of algebras manipulated by  classical Gr\"obner basis theory, such
as the Gr\"obner basis theory for solvable polynomial algebras or
their homomorphic images ([AL], [K-RW], [HT], [BGV], [Lev]), but may
hold a left Gr\"obner basis theory or a right Gr\"obner basis
theory. Inspired by such a fact,  we introduce more general skew
2-nomial algebras in this section, and, as the first step of having
a (one-sided or two-sided) Gr\"obner basis theory, we establish the
existence of a skew multiplicative $K$-basis for such algebras. As
to the existence of a (one-sided or two-sided) monomial ordering,
examples of Section 2 have told us that it has to be a matter of
examining concrete individual class of algebras.
\par

Throughout the current section $R$ denotes a $K$-algebra with a {\it
skew multiplicative $K$-basis} $\B$. Moreover, here and in what
follows, with respect to a skew multiplicative $K$-basis, a
one-sided or two-sided Gr\"obner basis theory means the one in the
sense of Definition 1.4. All notations used before are
maintained.\v5

We start by generalizing the sufficiency result of ([Gr2], Theorem
2.3).  Let $I$ be an arbitrary ideal of $R$. Then we may define a
relation $\sim_I$ on $\B\cup\{ 0\}$ as follows: if $v,u\in\B\cup\{
0\}$, then
$$v\sim_Iu\Longleftrightarrow \exists ~\lambda\in K^*~\hbox{such that}~v-\lambda u\in I.$$
{\parindent=0pt\par

{\bf 3.1. Lemma } $\sim$ is an equivalence relation on $\B\cup\{
0\}$.\vskip 6pt
{\bf Proof} Note that $K$ is a field, and $v\in I$ if and only if
$v\sim_I 0$. The verification of reflexivity, symmetry, and
transitivity for $\sim_I$ is straightforward. \QED}\v5

Recall that elements of $\B$ are called  monomials. If $h=v-\lambda
u\in R$ with $v,u\in\B$ and $\lambda\in K^*$, then we call $h$ a
{\it skew} 2-{\it nomial} in $R$; and thus we call the equivalence
relation $\sim_I$ obtained above a {\it skew} 2-{\it nomial
relation} on $\B\cup\{ 0\}$ determined by the ideal $I$.
{\parindent=0pt\v5

{\bf 3.2. Definition} An ideal $I$ of $R$ is said to be a {\it skew}
2-{\it nomial ideal} if $I$ is generated by monomials and skew
2-nomials. If $I$ is a skew 2-nomial ideal of $R$, then we call the
corresponding quotient algebra $A=R/I$ a {\it skew} 2-{\it nomial
algebra}. \v5

{\bf Remark} Instead of comparing to other references, for instance
[GI-1], [Laf], or [ES] in the commutative case, we used the phrase
``skew 2-nomial ideal" just for differing from the phrase ``2-nomial
ideal" used in [Gr2] (see a remark given after Theorem 3.4).\v5

{\bf 3.3. Lemma} Let $I$ be a skew 2-nomial ideal of $R$ and
$\sim_I$ the skew 2-nomial relation on $\B\cup\{ 0\}$ determined by
$I$. If $f=\sum_{j}\lambda_jw_j$ with $\lambda_j\in K^*$ and
$w_j\in\B$, then $f\in I$ if and only if for each equivalence class
$[w]$ of $\sim_I$, $\sum_{w_j\in [w]}\lambda_jw_j\in I$.\vskip 6pt
{\bf Proof} The sufficiency is clear. To get the necessity, suppose
$f=\sum_{j}\lambda_jw_j\in I$. Then since $\B$ is a skew
multiplicative $K$-basis and $I$ is an ideal, $f$ has a presentation
of the form
$$\begin{array}{rcl} f=\SUM^{}_{j}\lambda_jw_j&=&\SUM^{}_{p}\mu_p(v_p-\gamma_pu_p)
+\SUM^{}_{q}\eta_qs_q\\
&{~}&\hbox{with}~v_p-\gamma_pu_p,s_q\in I,~\mu_p,\gamma_p,\eta_q\in
K^*,\\
&{~}&\hbox{but}~v_p,u_p\not\in I~\hbox{for every}~p.\end{array}$$
Note that for $v$, $s\in\B$, if $s\in I$ but $v\not\in I$, then
$v\not\sim_I s$. So, comparing both sides of the above equality
subject to the equivalence relation $\sim_I$, if $[w]$ is any
equivalence class of $\sim_I$, then
$$\begin{array}{ll} \hbox{either}&\SUM^{}_{w_j\in [w]}\lambda_jw_j=
\sum_{v_p,u_p\in [w]}\mu_p(v_p-\gamma_pu_p)\in I,\\
\\
\hbox{or}&\SUM^{}_{w_j\in [w]}\lambda_jw_j=\sum_{s_q\in
[w]}\eta_qs_q\in I,\end{array}$$ as wanted.\QED}\v5

To make use of unified notation as before, if $I$ is an ideal of $R$
and $\pi$: $R\r R/I$ is the canonical surjection, then, for $f\in
R$, in what follows we write $\bar{f}$ for the canonical image of
$f$ in $R/I$, i.e., $\bar{f}=f+I$.{\parindent=0pt\v5

{\bf 3.4. Theorem} Let $I$ be a skew 2-nomial ideal of $R$, $A=R/I$,
and let $\sim_I$ be the skew 2-nomial relation on $\B\cup\{ 0\}$
determined by $I$.  Considering the image $\OV{\B}$ of $\B$ under
the canonical surjection $\pi$: $R\r A$, if we put
$\OV{\B}^*=\OV{\B}-\{ 0\}$, then the following statements hold.\par
(i) With respect to the inclusion order on subsets of $\OV{\B}^*$,
$\OV{\B}^*$ contains a maximal subset, denoted $\OV{\B}
^*_{\sim_I}$, with the property that
$$\OV{u_1},\OV{u_2}\in\OV{\B}^*_{\sim_I},~\OV{u_1}\ne \OV{u_2}
~\hbox{implies}~u_1\not\sim_Iu_2.$$
(ii) The maximal subset $\OV{\B}^*_{\sim_I}$ of $\OV{\B}^*$ obtained
in (i) above forms a skew multiplicative $K$-basis for the skew
2-nomial algebra $A$.\par
\vskip 6pt {\bf Proof} (i) The existence of $\OV{\B}^*_{\sim_I}$
with respect to the inclusion order on subsets of $\OV{\B}^*$
follows from Zorn's Lemma.\par
(ii) Since $\B$ is a skew multiplicative $K$-basis of $R$, if
$u,v\in\B$, then either $uv=0$ or $uv=\lambda w$ for some
$\lambda\in K^*$, $w\in\B$. Hence, for $\bar{u}$, $\bar{v}\in\OV{\B
}^*_{\sim_I}$, either $\bar{u}\bar{v}=0$, or $\bar{u}\bar{v}\ne 0$
and so $\bar{u}\bar{v}=\OV{uv}=\lambda\bar{w}$ with $\lambda\in K^*$
and $\bar{w}\in\OV{\B}^*$. If $\bar{w}\not\in\OV{\B }^*_{\sim_I}$,
then by the definition of $\OV{\B}^*_{\sim_I}$, there is some
$\bar{s}\in \OV{\B}^*_{\sim_I}$ such that $s\sim_Iw$. Consequently
$\bar{w}=\mu\bar{s}$ for some $\mu\in K^*$ and thereby
$\bar{u}\bar{v}=\lambda\mu \bar{s}$ with $\lambda\mu\in K^*$. Thus,
the skew multiplication property of $\OV{\B}^*_{\sim_I}$ is proved.
Furthermore, we proceed to show that $\OV{\B}^*_{\sim_I}$ forms a
$K$-basis of $A$. Again by the definition of $\OV{\B }^*_{\sim_I}$,
if $\bar{w}\in\OV{\B}^*$ and $\bar{w}\not\in\OV{\B }^*_{\sim_I}$,
then $\bar{w}=\lambda \bar{v}$ for some $\lambda\in K^*$ and
$\bar{v}\in\OV{\B}^*_{\sim_I}$. It follows that $\OV{\B
}^*_{\sim_I}$ spans $A$. To see that elements of $\OV{\B
}^*_{\sim_I}$ are linearly independent over $K$, suppose that
$\sum^{\ell}_{j=1}\lambda_j\OV{w_j}=0$ and that the $\OV{w_j}$ are
distinct, where $\lambda_j\in K$ and $\OV{w_j}\in\OV{\B
}^*_{\sim_I}$. Then $f=\sum^{\ell}_{j=1}\lambda_jw_j\in I$. By Lemma
3.3, $\sum_{w_j\in [u]}\lambda_jw_j\in I$ holds for every
equivalence class $[u]$ of $\sim_I$. But $w_j\in [u]$ implies
$w_j\sim_Iu$. Noticing that $\sim_I$ is an equivalence relation and
that all the $\OV{w_j}$ are distinct, it follows from the definition
of $\OV{\B}^*_{\sim_I}$ that different $w_j$'s are contained in
different equivalence classes. So $\lambda_jw_j\in I$ and
$\lambda_j\OV{w_j}=0$. Since $\OV{w_j}\ne 0$, it turns out that
$\lambda_j=0$ for all $j$. This completes the proof.\QED\v5

{\bf Remark} Let $R$ be a $K$-algebra with $K$-basis $\B$. Recall
from ([Gr2], Section 2) that an ideal $I$ of $R$ is said to be a
2-{\it nomial ideal} if $I$ can be generated by elements of the form
$b_1-b_2$ and $b$, where $b_1,b_2,b\in\B$. By ([Gr2], Theorem 2.3),
if $\B$ is a {\it multiplicative $K$-basis} of $R$ (i.e., $u,v\in\B$
implies $uv=0$ or $uv\in\B$) and $I$ is a 2-nomial ideal of $R$,
then, under the canonical epimorphism $\pi$: $R\r R/I$, the set $\OV
{\B}^*=\OV{\B}-\{ 0\}$ forms a multiplicative $K$-basis of $R/I$.
Since a multiplicative $K$-basis is trivially a skew multiplicative
$K$-basis and any 2-nomial ideal is trivially a skew 2-nomial ideal,
it is easy to see that in  the case that $\B$ is a multiplicative
$K$-basis  and $I$ is a 2-nomial ideal we must have $\OV{\B}^*=\OV
{\B}^*_{\sim_I}$, where the latter is the set obtained in Theorem
3.4(i) above. That is, actually the preceding Theorem 3.4
generalizes the sufficiency result of ([Gr2], Theorem 2.3).} \v5

Next, suppose further that there is a two-sided monomial ordering
$\prec$ on $\B$, that is, $(\B ,\prec)$ forms an admissible system
$(\B ,\prec )$ of $R$ and hence $R$ holds a two-sided Gr\"obner
basis theory (note that $\B$ is now a skew multiplicative $K$-basis
of $R$). Moreover, in this case if $I$ is an ideal of $R$, and
$A=R/I$, then it is easily seen that the set of normal monomials in
$\B$ (modulo $I$) is given by $N(I)=\B -\LM (I)$. Since the image
$\OV{N(I)}$ of $N(I)$ under the canonical algebra epimorphism $\pi$:
$R\r A$ forms a $K$-basis of $A$ and $\OV {N(I)}\subseteq
\OV{\B}^*$, by ([Gr2], Theorem 2.3), we first note the following
fact.{\parindent=0pt\v5

{\bf 3.5. Proposition} Let $I$ be a 2-nomial ideal of $R$ in the
sense of [Gr2]. Then $\OV{N(I)}=\OV{\B}^*$ and hence $\OV{N(I)}$ is
a multiplicative $K$-basis for $A=R/I$.\par\QED}\v5

Turning to the case of Theorem 3.4, we shall prove, under the
assumption that $R$ has an admissible system, that $\OV{N(I)}$ forms
a skew multiplicative $K$-basis for $A=R/I$.  {\parindent=0pt\v5

{\bf 3.6. Lemma} Suppose that $R$ has an admissible system $(\B
,\prec )$. Let $I$ be a skew 2-nomial ideal of $R$, $A=R/I$, and let
$N(I)$ be the set of  normal monomials in $\B$ (modulo $I$). With
notation as in Theorem 3.4,  the following statements hold.\vskip
6pt (i) For each $\bar{w}\in\OV{N(I)}$, there is a unique
$\lambda\in K^*$ and a unique $\bar{u}\in\OV{\B}^*_{\sim_I}$ such
that $\bar{w}=\lambda\bar{u}$; and if $\OV{w_1}$, $\OV{w_2}\in\OV
{N(I)}$ with $\OV{w_1}=\lambda_1\OV{u_1}$ and $\OV{w_2}=\lambda_2\pi
(u_2)$ respectively, then $\OV{w_1}\ne \OV{w_2}$ implies $\pi
(u_1)\ne \OV{u_2}$.\par
(ii) For each $\bar{u}\in \OV{\B}^*_{\sim_I}$, there is some $\bar
{w}\in\OV{N(I)}$ such that $\bar{u}=\mu\bar{w}$ with $\mu\in
K^*$.\vskip 6pt
{\bf Proof} To prove the assertions mentioned in (i) -- (ii), let us
bear in mind that $\OV{\B}^*_{\sim_I}$ is a $K$-basis of $A$ by
Theorem 3.4, and that $\OV{N(I)}$ is also a $K$-basis of $A$ by the
classical Gr\"obner basis theory.\par
(i) Note that $\OV{N(I)}\subseteq\OV{\B}^*$. Let $\bar{w}\in \OV
{N(I)}$. Then it follows from the definition of $\OV{\B
}^*_{\sim_I}$ that $\bar{w}=\lambda\bar{u}$ for some $\lambda\in
K^*$ and $\bar{u}\in\OV{\B}^*_{\sim_I}$. Since $\OV{\B}^*_{\sim_I}$
is a $K$-basis, the obtained expression is unique. If $\OV{w_1}, \OV
{w_2}\in\OV{N(I)}$ and $\OV{w_1}=\lambda_1\OV{u_1}$ and $\OV
{w_2}=\lambda_2\OV{u_2}$, then $\OV{u_1}=\OV{u_2}$ would clearly
imply the linear dependence of $\OV{w_1}$ and $\OV{w_2}$. Hence the
second assertion of (i) follows from the fact that $\OV{N(I)}$ is a
$K$-basis.  \par
(ii) Since $\OV{N(I)}$ is a $K$-basis for $A$, each $\bar{u}\in\OV
{\B}^*_{\sim_I}$ has a unique linear expression $\bar
{u}=\sum_{p=1}^m\mu_p\OV{w_p}$ with $\mu_p\in K^*$ and $\OV
{w_p}\in\OV{N(I)}$. By (i), $\bar{u}=\sum^m_{p=1}\mu_p\lambda_p\OV
{u_p}$ with $\lambda_p\in K^*$ and $\OV{u_p}\in\OV{\B }^*_{\sim_I}$,
in which $\OV{w_p}=\lambda_p\OV{u_p}$ and all the $\OV{u_p}$ are
distinct. If, without loss of generality, $\bar {u}=\OV{u_1}$, then
$\bar{u}=\lambda_1^{-1}\OV{w_1}$, as desired. \QED\v5

{\bf 3.7. Theorem} With notation and the assumption as in Theorem
3.4 and Lemma 3.6, $\OV{N(I)}$ is a skew multiplicative $K$-basis
for the skew 2-nomial algebra $A=R/I$.\vskip 6pt
{\bf Proof} If $\OV{w_1},\OV{w_2}\in\OV{N(I)}$ and $\OV{w_1}\OV{
w_2}\ne 0$, then by Lemma 3.6 and the fact that $\OV{\B
}^*_{\sim_I}$ is a skew multiplicative $K$-basis of $A$ (Theorem
3.4(ii)), there are $\lambda_1,\lambda_2,\lambda_3,\mu\in K^*$, $\OV
{u_1},\OV{u_2}\in\OV{\B}^*_{\sim_I}$ and $\OV{w_3}\in\OV{N(I)}$ such
that
$$\begin{array}{rcl} \OV{w_1}\cdot\OV{w_2}&=&\lambda_1\OV
{u_1}\lambda_2\OV{u_2}\\
&=&\lambda_1\lambda_2\lambda_3\OV
{u_3}\\
&=&\lambda_1\lambda_2\lambda_3\mu\OV{w_3}.\end{array}$$ It follows
that the $K$-basis $\OV{N(I)}$ of $A$ is a skew multiplicative
$K$-basis.\QED}\v5

As we will see below, indeed, a better proof of Theorem 3.7 may be
obtained by employing the division by a  Gr\"obner basis of $I$,
which involves the following interesting result (a noncommutative
analogue of ([ES], Proposition 1.1)). {\parindent=0pt\v5

{\bf 3.8. Proposition} Suppose that $R$ has an admissible system
$(\B ,\prec )$, and let $I$ be a skew 2-nomial ideal of $R$. Then
$I$ has a Gr\"obner basis $\G$ consisting of monomials and skew
2-nomials, i.e., elements of the form $v-\lambda u$, $s$ with
$\lambda\in K^*$ and $v,u,s\in\B$. \vskip 6pt
{\bf Proof} Let $G$ be any Gr\"obner basis of $I$. Then, since $I$
is generated by skew 2-nomials and monomials, as in the proof of
Lemma 3.3, each $g\in G$ has a presentation of the form
$$\begin{array}{rcl} g&=&\SUM^{}_{p}\mu_p(v_p-\gamma_pu_p)
+\SUM^{}_{q}\eta_qs_q\\
&{~}&\hbox{with}~v_p-\gamma_pu_p,s_q\in I,~\mu_p,\gamma_p,\eta_q\in
K^*,\\
&{~}&\hbox{but}~v_p,u_p\not\in I~\hbox{for every}~p.\end{array}$$
Let $g=\sum_{j}\lambda_jw_j$, where $\lambda_j\in K^*$ and
$w_j\in\B$. If $\LM (g)\not\in I$, then we claim that there is some
$w_j$ appearing in $\sum_{j}\lambda_jw_j$ such that  $w_j\ne\LM
(g)$, $\LM (g)-\mu w_j\in I$, and $\LM (\LM (g)-\mu w_j)=\LM (g)$,
where $\mu\in K^*$. Otherwise, considering the skew 2-nomial
relation $\sim_I$ on $\B\cup\{ 0\}$ determined by $I$, and letting
$[\LM (g)]$ be the equivalence class of of $\sim_I$ represented by
$\LM (g)$, by Lemma 3.3 we would have
$$\SUM^{}_{w_j\in [\LM (g)]}\lambda_jw_j=\LC (g)\LM (g)=
\sum_{v_p,u_p\in [\LM (g)]}\mu_p(v_p-\gamma_pu_p)\in I,$$
and consequently $\LM (g)\in I$, a contradiction. Thus, if we put
$$\G =\left\{ v-\lambda u,~s~\left |~\begin{array}{l}
v,u\in\B -I,~\lambda\in K^*,~v-\lambda u\in I,~s\in\B\cap I,~\hbox{such that}\\
\LM (v-\lambda u)=\LM (g_i),~ s=\LM (g_j)~\hbox{for some}~g_i,g_j\in
G\end{array}\right.\right \} ,$$
then $\LM (\G )=\LM (G)$ and hence $\G$ is a Gr\"obner basis for
$I$.\QED}\v5

For the convenience of statement, we call the Gr\"obner basis $\G$
obtained in Lemma 3.8 a {\it skew 2-nomial Gr\"obner basis} of the
skew 2-nomial ideal $I$.{\parindent=0pt\v5

{\bf 3.9. Lemma} Suppose that $R$ has an admissible system $(\B
,\prec )$. Let $I$ be a skew 2-nomial ideal of $R$ and $\G$ a skew
2-nomial Gr\"obner basis of $I$. If $u\in\B$ satisfying $u\not\in I$
and $u\not\in N(I)$, where $N(I)$ is the set of normal monomials in
$\B$ (modulo $I$), then $u$ has a presentation
$$\begin{array}{rcl} u&=&\displaystyle{\sum_{i,j}}\lambda_{ij}w_{ij}g_jv_{ij}+\lambda u'~\hbox{with}~\lambda_{ij},
\lambda\in K^*,~w_{ij},v_{ij}\in\B ,~g_j\in\G\\
&{~}&\hbox{and}~u'\in N(I)~\hbox{satisfying}~u'\prec u.\end{array}$$
{\bf Proof} It is simply a consequence of doing division by $\G$. As
$u\not\in I$ and $u\not\in N(I)$, there is some $g_i=u_i-\mu_is_i\in
\G$, such that $u=\lambda_iw_i\LM (g_i)v_i$ for some $\lambda_i\in
K^*$, $w_i,v_i\in\B$. Suppose $\LM (g_i)=u_i$. Then since $\B$ is a
skew multiplicative $K$-basis,
$$\begin{array}{rcl} u-\lambda_iw_ig_iv_i&=&\lambda_i\mu_iw_is_iv_i=\gamma_iu_1~\hbox{with}~\gamma_i\in K^*\\
&{~}&\hbox{and}~u_1\in\B~\hbox{satisfying}~u_1\prec u.\end{array}$$
Note that $u\not\in I$ implies $u_1\not\in I$. If $u_1\not\in N(I)$,
then repeat the division procedure for $u_1$. By the well-ordering
property of $\prec$, the proof is finished after a finite number of
reductions. \QED}\v5

Now, let $\pi$: $R\r A=R/I$ be the canonical epimorphism as before,
and suppose that $R$ has an admissible system $(\B ,\prec )$. If
$N(I)$ is the set of normal monomials in $\B$ (modulo $I$) and $\G$
is a skew 2-nomial Gr\"obner basis of $I$, then, since $\B$ is a
skew multiplicative $K$-basis, for any $u$, $v\in N(I)$, it follows
from Lemma 3.9 that
$$\begin{array}{l} \hbox{either}~\bar{u}\bar{v}=0, \\
\hbox{or}~\bar{u}\bar{v}=\OV{uv}=\lambda\bar
{w}~\hbox{with}~\lambda\in K^*,~w\in N(I).\end{array}$$ This shows
that $\OV{N(I)}$ is a skew multiplicative $K$-basis for the skew
2-nomial algebra $A$, that is, Theorem 3.7 is recaptured.\v5

Finally, concerning the existence of a Gr\"obner basis theory for a
skew 2-nomial algebra, let us mention the following
conclusion.{\parindent=0pt\v5

{\bf 3.10. Theorem} Let $R$ be a $K$-algebra with the skew
multiplicative $K$-basis $\B$, and let $I$ be a skew 2-nomial ideal
of $R$, $A=R/I$.\par (i) With notation as in Theorem 3.4, if there
is a (left, right, or two-sided) monomial ordering $\prec$ on
$\OV{\B}^*_{\sim_I}$, then the skew 2-nomial algebra $A$ has a
(left, right, or two-sided)  Gr\"obner basis theory.\par

(ii) With notation and the assumption as in Theorem 3.7, if there is
a (left, right, or two-sided) monomial ordering $\prec$ on
$\OV{N(I)}$, then the skew 2-nomial algebra $A$ has a (left, right,
or two-sided) Gr\"obner basis theory.}

\section*{4. Almost Skew 2-nomial Algebras}
Naturally, the results of Sections 2 -- 3, combined with the working
principle of [Li2] (as indicated before Proposition 2.6), motivate
us to introduce the class of algebras with the associated graded
algebra which is a skew 2-nomial algebra (such as the algebras given
in Example (8) of Section 2), namely the almost skew 2-nomial
algebras defined below. As in previous sections,  a one-sided or
two-sided Gr\"obner basis theory means the one in the sense of
Definition 1.4 or it is the one indicated before Theorem 2.10,
depending on whether a skew multiplicative $K$-basis is used or not.
All notations used before are maintained.\v5

Let $R=\oplus_{\gamma\in\Gamma}R_{\gamma}$ be a $\Gamma$-graded
$K$-algebra, where $\Gamma$ is an ordered semigroup by a total
ordering $<$, and we assume that $R$ has an admissible system $(\B
,\prec)$, where $\B$ is a {\it skew multiplicative $K$-basis} of $R$
consisting of $\Gamma$-{\it homogeneous elements}, and $\prec$ is a
two-sided monomial ordering on $\B$. Hence $R$ has a two-sided
Gr\"obner basis theory with respect to $(\B ,\prec)$.
\par Let $I$ be an ideal of $R$, $A=R/I$, and $A^{\Gamma}_{\rm
LH}=R/\langle\LH (I)\rangle$ the $\Gamma$-leading homogeneous
algebra of $A$ as defined in Section 2, where $\LH (I)$ is the set
of $\Gamma$-leading homogeneous elements of $I$. {\parindent=0pt\v5

{\bf 4.1. Definition}  With notation as above, if  $\langle\LH
(I)\rangle$ is a skew 2-nomial ideal, or equivalently, if
$A^{\Gamma}_{\rm LH}=R/\langle\LH (I)\rangle$ is a skew 2-nomial
algebra, then we call the quotient algbera $A=R/I$ an {\it almost
skew} 2-{\it nomial algebra}.}\v5

Before mentioning the main result of this section, let us
demonstrate, by referring to Examples (5) -- (7) of Section 2, how
to obtain an almost skew 2-nomial algebra in terms of Gr\"obner
bases in $R$. We begin by recording a property concerning the
homogeneous generators of the $\Gamma$-graded ideal $\langle\LH
(I)\rangle$, which is the same as  the property concerning the
monomial generators of the monomial ideal $\langle\LM (I)\rangle$ if
$\LH (~)$ once is replaced by $\LM (~)$.{\parindent=0pt\v5

{\bf 4.2. Lemma} With notation as above, let $h$ be a nonzero
$\Gamma$-homogeneous element of $R$. Then $h\in \langle\LH
(I)\rangle$ if and only if $h\in \LH (I)$.\vskip 6pt
{\bf Proof}  Since $h$ is a homogeneous element, if $h\in\langle\LH
(I)\rangle$, then
$$h=\sum_{i,j} H_{ij}\LH (f_i)T_{ij},~\hbox{where}~H_{ij},~T_{ij}~\hbox{are homogeneous elements and} ~f_i\in I.$$
If we put $f_i=\LH (f_i)+f_i'$, where deg$(f_i')<$ deg$(f_i)$, then
$f=\sum_{i,j} H_{ij}f_iT_{ij}\in I$ and
$$f=\sum_{ij} H_{ij}\LH (f_i)T_{ij}+\sum_{i,j} H_{ij}f_i'T_{ij}=h+\sum_{i,j} H_{ij}f_i'T_{ij}.$$
Hence it is clear that $h=\LH (f)\in\LH (I)$. The converse is
trivial.\QED \v5

{\bf 4.3. Proposition} Let $J$ be a $\Gamma$-graded ideal of $R$. If
$J$ is a skew 2-nomial ideal, then $J$ has a skew 2-nomial Gr\"obner
basis $\G$ (as obtained in Proposition 3.8) but consisting of
$\Gamma$-homogeneous elements.\vskip 6pt {\bf Proof} To be
convenient, let ${\rm H}^{\Gamma}(J)$ denote the set of
$\Gamma$-homogeneous elements of $J$. Noticing that $J$ is a
$\Gamma$-graded ideal, if $f=\sum_jr_{\gamma_j}\in R$ with
$r_{\gamma_j}\in R_{\gamma_j}$, then $f\in J$ if and only if
$r_{\gamma_j}\in J$ for all $j$.  Thus, since  $\B$ is a skew
multiplicative $K$-basis, exactly as in Proposition 1.3(i) we obtain
a homogeneous Gr\"obner basis for $J$ as follows
$$G=\left\{\left. h\in {\rm H}^{\Gamma}(J)~\right |~\hbox{if}~h'\in {\rm H}^{\Gamma}(J)
~\hbox{and}~h'\ne h,~\hbox{then}~\LM (h')\not |~\LM (h)\right\} .$$
Furthermore, noticing that $J$ is also a skew 2-nomial ideal,
actually as in the proof of Proposition 3.8, subject to $G$ above we
obtain a Gr\"obner basis of $J$ in the desired form
$$\G =\left\{ v-\lambda u,~s~\left |~\begin{array}{l}
v,u\in\B -J,~\lambda\in K^*,~v-\lambda u\in {\rm H}^{\Gamma}(J),~s\in\B\cap J,~\hbox{such that}\\
\LM (v-\lambda u)=\LM (h_i),~ s=\LM (h_j)~\hbox{for some}~h_i,h_j\in
G\end{array}\right.\right \} .$$ \par\QED\v5

{\bf 4.4. Corollary} Suppose that $R$ has a two-sided  admissible
system $(\B ,\prec_{gr})$ in which $\prec_{gr}$ is a $\Gamma$-graded
two-sided monomial ordering on $\B$, and let $A=R/I$ be an almost
skew 2-nomial algebra in the sense of Definition 4.1. Then $I$ has a
Gr\"obner basis $\mathscr{G}$ such that $\LH (\mathscr{G})$ is a
homogeneous skew 2-nomial Gr\"obner basis for the ideal $\langle\LH
(I)\rangle$ with respect to $(\B ,\prec_{gr})$. \vskip 6pt {\bf
Proof} Writing $J=\langle\LH (I)\rangle$, it follows from
Proposition 4.3 and Lemma 4.2 that $J$ has a skew 2-nomial Gr\"obner
basis $\G\subset\LH (I)$. Since $\prec_{gr}$ is a $\Gamma$-graded
monomial ordering on $\B$, by Proposition 2.4(ii), the desired
Gr\"obner basis $\mathscr{G}$ for $I$ is given by
$$\mathscr{G}=\{ f\in I~|~\LH (f)=g,~g\in\G\} .$$\par\QED}\v5

It follows from Corollary 4.4,  Theorem 2.10, Theorem 3.7, and
Theorem 3.10 that we can now write down the main result of this
section.{\parindent=0pt\v5

{\bf 4.5. Theorem} Suppose that $R$ has a two-sided admissible
system $(\B ,\prec_{gr})$ in which $\prec_{gr}$ is a $\Gamma$-graded
two-sided monomial ordering on $\B$. Let $I$ be an ideal of $R$,
$A=R/I$. If $I=\langle\mathscr{G}\rangle$ is generated by a
Gr\"obner basis $\mathscr{G}$ such that $\LH (\mathscr{G})$ consists
of monomials and $\Gamma$-homogeneous skew 2-nomials, then the
following statements hold.\par

(i) $\langle\LH (I)\rangle =\langle\LH (\mathscr{G})\rangle$ is a
$\Gamma$-graded skew 2-nomial ideal, the $\Gamma$-leading
homogeneous algebra $A^{\Gamma}_{\rm LH}=R/\langle\LH (I)\rangle$ of
$A$ is a skew 2-nomial algebra, and hence $A$ is an almost skew
2-nomial algebra.\par

(ii) Let $N(I)$ be the set of normal monomials in $\B$ (modulo $I$),
and write
$$\widehat{N(I)}=\{ \hat{u}=u+\langle\LH (I)\rangle~|~u\in N(I)\}$$
for the $K$-basis of $A^{\Gamma}_{\rm LH}$ (as in Section 2). Then
$\widehat{N(I)}$ is a skew multiplicative $K$-basis of
$A^{\Gamma}_{\rm LH}$. If furthermore there is a (left, right, or
two-sided) monomial ordering $\prec_{_{\rm LH(I)}}$ on
$\widehat{N(I)}$, then  $A^{\Gamma}_{\rm LH}$ has a (left, right, or
two-sided) Gr\"obner basis theory with respect to the admissible
system $(\widehat{N(I)},\prec_{_{\rm LH(I)}})$.\par

(iii) If $A^{\Gamma}_{\rm LH}$ is a domain and
$\prec_{\widehat{gr}}$ is a  $\Gamma$-graded (left, right, or
two-sided) monomial ordering on $\widehat{N(I)}$, then $A$ has a
(left, right, or two-sided)  Gr\"obner basis theory with respect to
the (left, right, or two-sided) admissible system $(\OV{N(I)},\prec
)$, where
$$\OV{N(I)}=\{ \OV u=u+I~|~u\in N(I)\}$$
is the $K$-basis of $A$ determined by $N(I)$ (as in Section 3), and
the monomial ordering $\prec$ on $\OV{N(I)}$ is defined subject to
the rule: for $\OV u$, $\OV v\in\OV{N(I)}$,
$\bar{u}\prec\bar{v}~\hbox{whenever}~\hat{u}\prec_{\widehat{gr}}\hat{v}$.}
\v5

\section*{5. Some Open Problems}
As one may see from previous sections,  the (left, right, or
two-sided) Gr\"obner basis theory introduced for algebras with a
skew multiplicative $K$-basis not only is in principle consistent
with the well-known existed Gr\"obner basis theory, but also covers
certain new classes of algebras (such as two subclasses of quantum
binomial algebras, down-up algebras, and algebras with the
associated graded algebra of such types) which are beyond the scope
of algebras treated by the well-known existed Gr\"obner basis
theory. However, the already obtained results concerning (almost)
skew 2-nomial algebras are obviously far from being complete, for
instance, even if for an algebra of the type as described in
Examples (1), respectively Example (2) of Section 2 ( which is a
Noetherian quantum binomial algebra whenever $\lambda,\mu ,
\gamma\in K^*$ and $\gamma^2=1$, respectively  $\lambda^2=1$), we do
not know if there is an implementable analogue of Buchberger
Algorithm for producing a finite (left or right) Gr\"obner basis;
and moreover, proper subclasses of practical (almost) skew 2-nomial
algebras, that are not the types we have considered in Section 2 but
have a (left, right, or two-sided) monomial ordering, need to be dug
further. Therefore, we finish this paper by summarizing several
 open problems related to preceding sections: {\parindent=0pt\par

{\bf 1.} Under the Noetherian assumption, is there an analogue of
Buchberger Algorithm for the algebras described in Corollary 2.5(i)
-- (iii) respectively? \par

{\bf 2.} For the algebras described in Corollary 2.5(i) -- (iii),
especially  for those described in (i) -- (ii) (including many
quantum binomial algebras), is there a computational module theory
(representation theory) in terms of (left, right, or two-sided)
Gr\"obner bases as that developed in [Gr1] and [Gr2]?\par

{\bf 3.} Find more subclasses  of skew 2-nomial algebras, in which
each algebra is not the type as described in Theorem 2.2 and Theorem
2.3 respectively, but has a (left, right, or two-sided) monomial
ordering. \par

{\bf 4.} Find more subclasses  of almost skew 2-nomial algebras, in
which the $\Gamma$-leading homogeneous algebra $A^{\Gamma}_{\rm LH}$
of each algebra $A$ is not the type as described in Theorem 2.7 and
Theorem 2.8 respectively, but satisfies the conditions of Theorem
2.10.\par

{\bf 5.} For algebras as described in Example (8) of Section 2, if
the Gr\"obner basis theory of $A^{\NZ}_{\rm LH}$ is algorithmically
realizable, is it possible to realize the lifted Gr\"obner basis
theory of $A$ in an algorithmic way?\par

{\bf 6.} In view of Remark (2) given before Theorem 2.10, is there
an algorithmically realizable Gr\"obner basis theory for the down-up
algebras of the type $A(\alpha ,0,\gamma )$?\par

{\bf 7.} Let $\G$ be a Gr\"obner basis of the ideal
$I=\langle\G\rangle$ in the free $K$-algebra $\KS =K\langle
X_1,...,X_n\rangle$ with respect to some $\NZ$-graded monomial
ordering, such that $\LH (\G )$ is of the type as described in
Corollary 2.5(i) or Corollary 2.5(ii), but with $\Omega =\emptyset$
and all the $\lambda_{ji}\ne 0$,  and such that the $\NZ$-leading
homogeneous algebra $A^{\NZ}_{\rm LH}=\KS /\langle\LH (\G )\rangle$
of $A=\KS /I$ is a quantum binomial algebra in the sense of [Laf]
and [G-I2]. By referring to [Laf], [G-I2] and [Li2], explore the
structural properties of $A$ via $A^{\NZ}_{\rm LH}$ in a
computational way.

\v5

 \centerline{References}
\parindent=1.3truecm

\re{[AL]} J. Apel and W. Lassner, An extension of Buchberger's
algorithm and calculations in enveloping fields of Lie algebras,
{\it J. Symbolic Comput}., 6(1988), 361--370.
\item{[Ap]} J. Apel, Computational ideal theory in finitely generated extension rings,
{\it Theoretical Computer Science},  244(2000) 1 -- 33.
\item{[Ben]} G.~Benkart, Down-up algebras and Witten's deformations
of the universal enveloping algebra of $sl_2$, {\it Contemp. Math.},
224(1999), 29--45.
\item{[Bu]} B. Buchberger, Gr\"obner bases: An algorithmic method in
polynomial ideal theory, in: {\it Multidimensional Systems Theory}
(N.K. Bose, ed.), Reidel Dordrecht, 1985, 184--232.
\item{[BGV]} J. Bueso, J. G\'omez--Torrecillas, and A. Verschoren,
{\it Algorithmic methods in non-commutative algebra}: {\it
Applications to quantum groups}. Kluwer Academic Publishers, 2003.
\item{[CM]} P.~Carvalho and M.~Musson, Down-up algebras and their
representation theory, {\it J. Alg.}, 228(2000), 286--310.
\item{[CS]} T. Cassidy and B. Shelton, Basic properties of
generalized down-up algebras, {\it J. Alg}., 279(2004), 402-421
\item{[ES]} D. Eisenbud and B. Sturmfels, Binomial ideals,
{\it Duke Math. J.}, 1(84)(1996), 1 -- 45.
\item{[FFG]} D.R. Farkas, C. Feustel and E.L. Green, Synergy in the theories of Gr\"obner bases
and path algebras, {\it Canad. J. Math.}, 45(1993), 727 -- 739.
\item{[G-I1]} T. Gateva-Ivanova, Noetherian properties of skew polynomial rings with binomial
relations, {\it Trans. Amer, Math. Sco.}, 1(343)(1994), 203 -- 219
\item{[G-I2]} T. Gateva-Ivanova, Binomial skew polynomial rings, Artin-Schelter
regularity, and binomial solutions of the Yang-Baxter equation,
arXiv:0909.4707.

\item{[GL]} K.R. Goodearl and E.S. Letzter, Prime and primitive spectra
of multiparameter quantum affine spaces, in: Trends in Ring Theory
(MisKolc, 1996) (V. Dlab and L. marki, eds.), {\it Canad. Math. Soc.
Conf. Proc. Series} 22(1998), 39 -- 58.

\item{[Gr1]} E.L. Green, Noncommutative Gr\"obner bases and projective resolutions, in:
{\it Proceedings of the Euroconference Computational Methods for
Representations of Groups and Algebras, Essen, 1997}, (Michler,
Schneider, eds), Progress in Mathematics, Vol. 173, Basel,
Birkha\"user Verlag, 1999, 29 -- 60.
\item{[Gr2]} E.L. Green, Multiplicative bases, Gr\"obner bases, and right Gr\"obner bases,
{\it J. Symbolic Comput.}, 29(2000), 601 -- 623.
\item{[HT]} D. Hartley and P. Tuckey, Gr\"obner Bases in Clifford and
Grassmann Algebras, {\it J. Symb. Comput.}, 20(1995), 197--205.
\item{[KMP]} E.~Kirkman, I.M.~Musson, and D.S.~Passman, Noetherian
down-up algebras, {\it Proc. Amer. Math. Soc.}, 127(1999),
2821--2827.
\item{[K-RW]} A.~Kandri-Rody and V.~Weispfenning, Non-commutative
Gr\"obner bases in algebras of solvable type, {\it J. Symbolic
Comput.}, 9(1990), 1--26.
\item{[Laf]} G. Laffaille, Quantum binomial algebras, Colloquium on Homology and
Representation Theory (Spanish) (Vaquerfas, 1998). Bol. Acad. Nac.
Cienc. (C\'ordoba) 65 (2000), 177--182.
\item{[Le-B]} L.~Le Bruyn, Conformal $sl_2$ enveloping algebras, {\it
Comm. Alg.}, 23(1995), 1325--1362.
\item{[Lev]} V. Levandovskyy, {\it Non-commutative Computer Algebra for
Polynomial Algebra}: {\it Gr\"obner Bases, Applications and
Implementation}, Ph.D. Thesis, TU Kaiserslautern, 2005.
\item{[Li1]} H. Li, {\it Noncommutative Gr\"obner Bases and
Filtered-Graded Transfer}, LNM, 1795, Springer-Verlag, 2002.
\item{[Li2]} H. Li, $\Gamma$-leading homogeneous algebras and Gr\"obner bases,
in: {\it Recent Developments in Algebra and Related Areas}, Advanced
Lectures in Mathematics, Vol. 8, International Press, Boston, 2009,
155 -- 200.\newpage
\item{[Man]} Yu.I.~Manin, {\it Quantum Groups and Noncommutative
Geometry}, Les Publ. du Centre de R\'echerches Math., Universite de
Montreal, 1988.
\item{[Mor]} T. Mora, An introduction to commutative and noncommutative
Gr\"obner Bases, {\it Theoretic Computer Science}, 134(1994),
131--173.

\end{document}